\documentclass[a4paper]{article}
\usepackage[margin=1.05in]{geometry}

\usepackage{authblk}
\usepackage{setspace}
\usepackage{enumitem}

\usepackage{array}
\usepackage{makecell}
\usepackage{amsmath, amssymb, booktabs, mathtools}
\usepackage{subfig}
\usepackage[font=small,labelfont=bf]{caption}
\usepackage{threeparttable}
\usepackage{float}
\usepackage[square, sort&compress, numbers]{natbib}
\usepackage{bm}
\usepackage{stmaryrd}
\usepackage{multirow}

\usepackage{graphicx} 
\usepackage{hyperref} 
\usepackage[nameinlink,capitalize]{cleveref}
\hypersetup{
  colorlinks   = true,
  urlcolor     = blue,
  linkcolor    = blue,
  citecolor    = blue
}
\newcommand{\correspondingA}{*}

\title{Mode veering and symmetry-protected crossings in conservative elastic waveguides: unified perturbation-theoretic interpretation and adaptive tracking}

\author[]{Dong Xiao\textsuperscript{*}\href{https://orcid.org/0009-0006-7609-7832}{\includegraphics[scale=0.04]{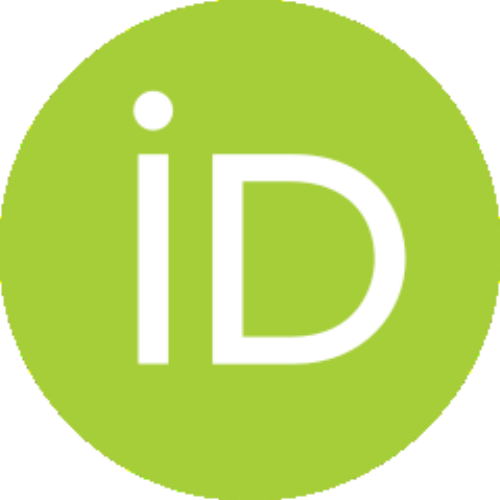}}}
\author[]{Zahra Sharif-Khodaei\href{https://orcid.org/0000-0001-5106-2197}{\includegraphics[scale=0.04]{orcidicon.pdf}}}
\author[]{M. H. Aliabadi\href{https://orcid.org/0000-0002-2883-2461}{\includegraphics[scale=0.04]{orcidicon.pdf}} }

\affil[]{\normalsize \slshape  Department of Aeronautics, Imperial College London, South Kensington, London SW7 2AZ, United Kingdom.}
\date{\vspace*{-2\baselineskip}}  
\date{\vspace*{-1.5cm}}

\begin{document}
\maketitle
\footnotetext[0]{ \textsuperscript{\correspondingA}Corresponding author}
\footnotetext[1]{Email addresses: d.xiao21@imperial.ac.uk (D. Xiao); z.sharif-khodaei@imperial.ac.uk (Z. Sharif-Khodaei); m.h.aliabadi@imperial.ac.uk (M.H. Aliabadi)}


\renewcommand{\abstractname}{Abstract}
\begin{abstract}{\normalsize \onehalfspacing
Accurate mode tracking is essential for elastic waveguide dispersion analysis in ultrasonic nondestructive evaluation and structural health monitoring. Its reliability, however, deteriorates near mode veering and closely spaced eigenvalues, where rapid eigenvector exchange causes mode misidentification. Although widely observed, the quantitative relationship between veering, eigenvector evolution, and tracking robustness has not been systematically established.  By specializing classical perturbation theory to the single‑parametric Hermitian SAFE eigenproblem---exemplified by conservative elastic waveguides---we obtain explicit expressions for eigenvector derivatives and modal coupling strength.  This yields a unified, quantitative interpretation of mode veering, symmetry-protected crossings, and degeneracies, and clarifies their distinct tracking implications: eigenvector sensitivity scales inversely with the eigengap, explaining modal repulsion and the degradation of correlation‑based tracking near avoided crossings, whereas symmetry-protected crossings remain benign because symmetry‑induced decoupling preserves smooth eigenvector evolution, and symmetry-protected degeneracies require rotation‑invariant subspace tracking.  A numerical consistency condition and an existence result for a critical step size are then derived, motivating a two‑level adaptive strategy with an a posteriori error indicator that separates numerical tracking consistency from symmetry‑based physical correctness. Numerical examples validate the theoretical predictions and demonstrate improved robustness in regions of strong modal interaction, providing practical guidance for reliable dispersion calculations and ultrasonic inspection.
}
\end{abstract}
\renewcommand{\abstractname}{Keywords}
\begin{abstract}{Guided wave dispersion; Mode veering; Symmetry-protected crossing; Perturbation theory; Adaptive step-size strategy; Non-destructive testing}
\end{abstract}
\newpage

\section{Introduction}

Guided waves in bounded elastic structures exhibit strong dispersion and
multimodal characteristics, making the accurate construction of dispersion
curves fundamental to ultrasonic nondestructive evaluation (NDE), structural
health monitoring (SHM), and material characterization~\cite{rose_ultrasonic_2014,
giurgiutiu_structural_2014}.  Reliable mode identification is a prerequisite for
transducer design, signal interpretation, defect imaging, and inverse problems
across applications ranging from aerospace composites to civil infrastructure
and biomedical imaging~\cite{aliabadi_structural_2018, xiao_hybrid_2024,
kaczmarek_noncontact_2017, laugier_bone_2011}.  However, the physical
interaction of dispersion branches poses a persistent challenge.  In
\emph{mode veering} (also known as avoided crossing) regions, eigenvalues
approach each other without intersecting while eigenvectors undergo rapid
exchange over a narrow parameter interval, causing mode misidentification that
compromises inspection reliability.  In contrast, when the waveguide possesses
reflection symmetry, modes belonging to different symmetry families may exhibit
true eigenvalue crossings without eigenvector mixing---a behaviour that, as
will be shown, reflects the underlying symmetry-induced decoupling of the
governing operator.  Understanding how these two types of modal interaction
influence eigenvector evolution and tracking reliability is essential for both
wave physics and practical inspection techniques.

The phenomenon of mode veering has been recognized since the seminal work of
Perkins and Mote~\cite{perkins_comments_1986} as a generic feature of coupled
eigenvalue problems.  In elastic waveguides, Zharnikov and
Syresin~\cite{zharnikov_repulsion_2015} demonstrated that quasidipole modes
repel when geometric symmetry is broken and undergo true crossings when the
symmetry is exact, attributing this to the perturbation‑theoretic splitting of
degenerate eigenvalues.  In the present work, drawing on the Wigner–von Neumann
non‑crossing rule~\cite{von_neumann_no_1929, hatton_noncrossing_1976,
mead_noncrossing_1979}, we adopt the more precise designation
\emph{symmetry-protected crossing} to emphasize that such crossings are
protected by symmetry: the governing operator is block-diagonalized, forcing
the inter-modal coupling to vanish identically and preventing the crossing from
being opened into a veering by any symmetry-preserving perturbation.  Despite
these insights, a quantitative relationship between the eigenvalue gap and the
eigenvector exchange rate---essential for predicting tracking reliability---has
not been systematically established.  Although both phenomena have been
extensively reported, the conditions under which symmetry-protected crossings
remain numerically benign, and the gap‑sensitivity relationship that governs
veering, have lacked a unified theoretical treatment.

Part of the difficulty stems from early dispersion formulations.  Systematic
approaches---including the global matrix method (GMM)~\cite{lowe_matrix_1995},
transfer matrix methods (TMM)~\cite{maghsoodi_calculation_2014}, and stiffness
matrix methods (SMM)~\cite{rokhlin_stable_2002}---formulate the waveguide
problem by searching for roots of $\det \mathbf{D}(\omega, \xi) = 0$, often
without simultaneously constructing or tracking the eigenvectors.  Huber and
Sause~\cite{huber_classification_2018, huber_classification_2023} improved this
practice by partitioning roots into modal families based on laminate symmetry,
effectively suppressing mode jumping across families.  However, this strategy
is restricted to plate‑like structures with known lay‑up symmetries, collapses
in asymmetric configurations where pure families no longer exist, and provides
no mechanism to resolve veering within a single family.  The fundamental
limitation is that the root‑finding pipeline operates on the scalar dispersion
relation, so eigenvector evolution is not directly available to constrain the
tracking.

A paradigmatic shift occurred with the Semi‑Analytical Finite Element (SAFE)
method~\cite{bartoli_modeling_2006, marzani_semi-analytical_2008} and spectral
collocation techniques~\cite{quintanilla_guided_2015, adamou_spectral_2004},
which project the wave operator onto a finite‑dimensional basis, yielding the
parameterized generalized Hermitian eigenproblem
\begin{equation}
[\mathbf{K}(k) - \omega^2 \mathbf{M}] \mathbf{q} = \mathbf{0},
\label{eq: safe_eigenproblem}
\end{equation}
where $\mathbf{K}(k) = \mathbf{K}_1 + i k \mathbf{K}_2 + k^2 \mathbf{K}_3$ is
Hermitian for lossless media and any real wavenumber $k$. The complete eigenbasis is computed at each wavenumber, providing explicit, continuous eigenvectors that make
perturbation‑theoretic analysis directly applicable.  SAFE extends naturally
to arbitrary cross‑sections~\cite{duan_investigation_2019, hayashi_guided_2003,
nguyen_ultrasonic_2013} and has been applied to piezoelectric~\cite{cortes_dispersion_2008},
prestressed~\cite{quiroga_estimation_2017, mazzotti_guided_2012},
fluid‑coupled~\cite{castaings_finite_2008}, and functionally graded~\cite{nanda_semi-analytical_2025}
waveguides. Implementations range from root‑finding‑based tools such as
DISPERSE~\cite{pavlakovic_disperse_1997} and the Dispersion
Calculator~\cite{huber_dispersion_2024, huber_stiffness_2024} to SAFE‑based
solvers such as SAFEDC~\cite{liu_modelling_2024} and
GUIGUW~\cite{bocchini_graphical_2011}; comparative studies~\cite{orta_comparative_2022,
quiroga_evaluation_2025, barazanchy_comparative_2017} have confirmed their
accuracy.  Despite these advances, robust mode tracking remains an unresolved
challenge, particularly in regions of strong modal interaction.

Within spectral formulations, constructing continuous dispersion curves
requires linking eigenpairs across parameter steps---mode tracking.  Two
strategies dominate.  Similarity‑based methods, principally the Modal
Assurance Criterion (MAC)~\cite{allemang_correlation_1982}, match eigenvectors
by correlation~\cite{marwala_finite-element-model_2010}, but their reliability
degrades in veering regions where rapid eigenvector rotation exceeds the
resolving capacity of a fixed step~\cite{leissa_curve_1974}.  Continuation
methods~\cite{maruyama_continuation_2025} trace each branch individually with
predictor–corrector schemes, but are sequential, require initial guesses, and
become ill‑conditioned near degeneracies.  Crucially, both paradigms treat mode
tracking as a numerical matching problem rather than a physical dynamics
problem: neither relates tracking failure to the underlying eigenvalue gap, nor
exploits the Hermitian structure to establish theoretical guarantees.

Recent work has begun to bridge this gap.  Gravenkamp et
al.~\cite{gravenkamp_notes_2023} exploited the Hermitian structure of SAFE to
separate modal families through block‑diagonalisation, eliminating cross‑family
coupling and reducing tracking to independent subproblems.  Yet their analysis
stops at classification: it does not explain the dynamics \emph{within}
families---why eigenvectors exchange rapidly in veering regions, how the
exchange rate relates to the eigenvalue gap, or whether tracking consistency
can be guaranteed.  The step from structural classification to physical
dynamics, and from empirical refinement to theoretically grounded adaptive
sampling, remains unmade.

This work addresses this gap by returning to the Hermitian structure inherent
in the SAFE formulation for lossless, traction‑free waveguides.  Because the
stiffness matrix $\mathbf{K}(k)$ is Hermitian for all real $k$, classical
perturbation theory~\cite{kato_perturbation_1995} can be specialised to yield
an explicit expression for the eigenvector derivative. We show that the rate
of eigenvector exchange scales inversely with the eigenvalue gap.  This yields
a unified interpretation of mode veering as a predictable consequence of local
spectral properties, and explains why correlation-based tracking degrades
precisely where the eigengap becomes small.  The same perspective shows
symmetry-protected crossings to be numerically benign: symmetry-induced
decoupling eliminates coupling terms, guaranteeing smooth eigenvector evolution
through the crossing.  Together, mode veering and symmetry-protected crossings
embody the two generic outcomes of the Wigner–von Neumann non‑crossing rule for
single‑parametric Hermitian eigenvalue problems.  These results directly
motivate a practical, physics-informed adaptive sampling strategy driven by a
numerical consistency condition (NCC) and an a posteriori error indicator,
together with a rotation-invariant subspace measure that handles
symmetry-protected degeneracy in axisymmetric waveguides.

The specific contributions are:
\begin{enumerate}
    \item \textbf{Specialization of eigenvector perturbation theory to SAFE.}
    We establish the inverse eigengap–sensitivity relationship, providing a
    unified, quantitative interpretation of mode veering that directly links
    the eigengap to tracking reliability.
    \item \textbf{Unified interpretation of symmetry-protected crossings.}
    Within the same framework, we show that symmetry‑protected crossings
    correspond to vanishing modal coupling, guaranteeing bounded eigenvector
    derivatives.  For continuous symmetries producing degenerate subspaces, we
    introduce subspace MAC as a rotation‑invariant similarity measure.
    \item \textbf{Existence of a critical step size for MAC‑based tracking.}
    We show that, for any veering region with a positive eigengap, there
    exists a critical step size $\Delta k_{\max}(k)$ below which the NCC holds.
    This provides a direct, existence-based theoretical justification for
    adaptive step‑size control.
    \item \textbf{Two-level adaptive refinement strategy.} We develop an adaptive
    algorithm that separates numerical tracking consistency (Level~1,
    NCC-driven) from symmetry-based physical correctness (Level~2,
    symmetry-penalized).  A unified framework handles both non‑degenerate modes
    and degenerate subspaces with a built‑in error indicator.
    \item \textbf{Comprehensive numerical validation.} The method is validated
    on four examples---symmetric and unsymmetric laminated plates, an L‑shaped
    bar, and a steel pipe---with comparisons against SAFEDC and the Dispersion
    Calculator.
\end{enumerate}

The paper is organized as follows. \cref{sec:methodology} presents the SAFE
formulation, MAC-based tracking, perturbation analysis, classification of modal
interactions, symmetry protection interpretation, existence of a critical step
size (Proposition~1), and adaptive algorithm. \cref{sec:numerical} provides
numerical validations, discusses computational efficiency, compares the
fundamental distinctions between available methods and analyzes waveguides
beyond exact symmetry. \cref{sec:conclusion} summarizes the main conclusions,
acknowledges limitations, and outlines directions for future work.

\section{Methodology}
\label{sec:methodology}
\subsection{Semi-Analytical Finite Element (SAFE)}

The Semi-Analytical Finite Element (SAFE) method exploits the translational
symmetry of a uniform waveguide to reduce the three-dimensional wave propagation
problem to a two-dimensional finite element discretization of the cross-section
\cite{bartoli_modeling_2006}.  For a waveguide extending in the $x$-direction
with arbitrary cross-section in the $(y,z)$-plane, the displacement field within
an element is expressed as
\begin{equation}
    \mathbf{u}^{(e)}(x,y,z,t) = \mathbf{N}^{(e)}(y,z)\,\mathbf{q}^{(e)}\,
                               e^{i(kx - \omega t)},
\end{equation}
where $\mathbf{N}^{(e)}$ contains the shape functions, $\mathbf{q}^{(e)}$ the
element nodal displacements, $k$ the wavenumber, $\omega$ the angular frequency,
and $t$ time.

Substituting this ansatz into the elastodynamic equations and applying the
principle of virtual work yields a quadratic eigenvalue problem in $k$:
\begin{equation}
    (\mathbf{K}_1 + ik\mathbf{K}_2 + k^2\mathbf{K}_3 - \omega^2\mathbf{M})\mathbf{q}
    = \mathbf{0},
\end{equation}
where $\mathbf{q}$ collects the global cross-sectional degrees of freedom.  The
matrices are assembled from element contributions
\begin{equation}
    \begin{aligned}
        \mathbf{K}_1 &= \int_A \mathbf{B}_0^{T}\mathbf{C}\mathbf{B}_0\,dA, \quad
        \mathbf{K}_2 = \int_A (\mathbf{B}_0^{T}\mathbf{C}\mathbf{B}_1 -
                              \mathbf{B}_1^{T}\mathbf{C}\mathbf{B}_0)\,dA, \\
        \mathbf{K}_3 &= \int_A \mathbf{B}_1^{T}\mathbf{C}\mathbf{B}_1\,dA, \quad
        \mathbf{M}   = \int_A \rho\mathbf{N}^{T}\mathbf{N}\,dA,
    \end{aligned}
\end{equation}
with strain–displacement matrices $\mathbf{B}_0$, $\mathbf{B}_1$, real material
stiffness $\mathbf{C}$ (lossless media), shape functions $\mathbf{N}$, mass
density $\rho$, and cross-sectional area $A$.  $\mathbf{K}_1$, $\mathbf{K}_3$,
and $\mathbf{M}$ represent cross-sectional, axial, and inertial contributions,
respectively, while $\mathbf{K}_2$ couples axial and cross-sectional
deformations.  For lossless materials with traction‑free boundaries,
$\mathbf{K}_1$, $\mathbf{K}_3$, $\mathbf{M}$ are symmetric and $\mathbf{K}_2$ is
skew‑symmetric.

We adopt wavenumber sweeping: $k$ is prescribed as a real parameter and the
corresponding eigenvalues $\lambda = \omega^2$ are obtained from the generalized
Hermitian eigenproblem
\begin{equation}\label{eq:eigenproblem}
    \mathbf{K}(k)\,\mathbf{q} = \lambda\,\mathbf{M}\,\mathbf{q},
    \qquad \mathbf{K}(k) = \mathbf{K}_1 + i k\mathbf{K}_2 + k^2\mathbf{K}_3 .
\end{equation}
For real $k$, $\mathbf{K}(k)$ is complex Hermitian while $\mathbf{M}$ is real,
symmetric, positive definite, and independent of $k$; this Hermitian structure
reflects energy conservation in the lossless waveguide.  All eigenvalues are
therefore real and the eigenvectors form a complete $\mathbf{M}$-orthonormal
basis.  The principal difficulty is not computing the eigenpairs at each $k$, but
tracking them continuously so that every dispersion branch follows the physical
evolution of a single mode.

\subsection{Mode tracking via the Modal Assurance Criterion (MAC)}
\label{sec:tracking}

\subsubsection{Definition of MAC}

For $\mathbf{M}$-orthonormal eigenvectors ($\mathbf{q}_i^H\mathbf{M}\mathbf{q}_j = \delta_{ij}$), the Modal Assurance Criterion (MAC)~\cite{pastor_modal_2012} between $\mathbf{q}_i(k)$ and $\mathbf{q}_j(k+\Delta k)$ is defined as
\begin{equation}\label{eq:MAC}
    \mathrm{MAC}[\mathbf{q}_i(k), \mathbf{q}_j(k+\Delta k)] =
    \bigl| \mathbf{q}_i^{H}(k)\,\mathbf{M}\,\mathbf{q}_j(k+\Delta k) \bigr|^2,
\end{equation}
taking values in $[0,1]$.  A MAC close to $1$ indicates strong physical similarity of the modal displacement fields; a decay of the self-MAC across adjacent wavenumbers therefore signals a rapid exchange of modal character, as occurs in veering.

\subsubsection{Numerical consistency condition (NCC)}
\label{sec:NCC}

To establish correspondences between the mode sets at $k$ and $k+\Delta k$, we build a cost matrix $\mathbf{C}$ with entries $C_{ij}=1-\mathrm{MAC}[\mathbf{q}_i(k),\mathbf{q}_j(k+\Delta k)]$ and apply the Hungarian algorithm~\cite{kuhn_hungarian_1955} to obtain the globally optimal one-to-one matching $\sigma$ that minimizes $\sum_i C_{i,\sigma(i)}$.

A sufficient condition for this assignment to preserve label continuity --- i.e., for the algorithm to return the identity permutation $\sigma(i)=i$ without branch jumping --- is that the self-MAC strictly exceeds every cross-MAC both row-wise and column-wise:
\begin{equation}\label{eq:NCC}
\begin{aligned}
    \mathrm{MAC}[\mathbf{q}_i(k), \mathbf{q}_i(k+\Delta k)] &> \mathrm{MAC}[\mathbf{q}_i(k), \mathbf{q}_j(k+\Delta k)] \quad \forall j\neq i,\\
    \mathrm{MAC}[\mathbf{q}_i(k), \mathbf{q}_i(k+\Delta k)] &> \mathrm{MAC}[\mathbf{q}_j(k), \mathbf{q}_i(k+\Delta k)] \quad \forall j\neq i.
\end{aligned}
\end{equation}
When \cref{eq:NCC} holds, every diagonal cost $C_{ii}$ is strictly smaller than all other entries in its row and column, forcing the identity permutation to be the unique optimal assignment.  The conditions are sufficient but not necessary; even when violated, the Hungarian algorithm may still recover the correct matching if the global cost structure favours the identity.  However, the row condition alone does not guarantee $\sigma(i)=i$ (especially when $M>N$), because a permutation could satisfy row minima while violating column minima.  Both conditions are therefore required for a rigorous sufficient criterion at the algorithmic level.

We refer to \cref{eq:NCC} as the \emph{Numerical Consistency Condition} (NCC).  It guarantees that the Hungarian algorithm returns the identity permutation and that mode labels remain continuous across adjacent steps.  Crucially, the NCC is a purely numerical criterion: it ensures label continuity, not physical mode identity.  In a severe veering zone, the condition can be satisfied even while the physical modes have effectively exchanged character, provided the eigenvalue gap remains above numerical noise.  Physical correctness in ambiguous regions must therefore be established by one of two complementary mechanisms: (i)~symmetry-based constraints (e.g., the Wigner--von Neumann non-crossing rule), which forbid crossings between modes of the same symmetry class; or (ii)~adaptive numerical resolution, which refines the wavenumber step until the local eigenvector variation is fully captured by the MAC (see Proposition~1 and \cref{sec:existence}).  The latter directly motivates the adaptive strategy of \cref{sec:adaptive}.

\subsection{Perturbation analysis of eigenvectors}
\label{sec:perturbation}

To understand how eigenvectors evolve with wavenumber $k$ — and ultimately to
prove the existence of a step size that guarantees numerically consistent tracking — we now
specialise classical perturbation theory to the SAFE eigenproblem.  The analysis
is based on two standard assumptions:
\begin{enumerate}[label=(\roman*)]
    \item $\mathbf{K}(k)$ is analytic on $[k_{\min},k_{\max}] \subseteq \mathbb{R}$;
    \item the eigenvalues are simple: $\lambda_i(k) \neq \lambda_j(k)$ for all
          $i\neq j$ and all $k$ in this interval.
\end{enumerate}
Under these conditions the eigenvectors of $\mathbf{K}\mathbf{q}_i =
\lambda_i\mathbf{M}\mathbf{q}_i$ can be chosen to vary analytically with $k$ and are
unique up to a phase factor \cite{kato_perturbation_1995}.  (Although Kato's
treatment assumes a standard inner product, the $\mathbf{M}$-weighted orthogonality
used here is handled by first transforming to a standard eigenproblem via
Cholesky factorization of $\mathbf{M}$.) 

Differentiating the eigenproblem with respect to $k$ and projecting onto another
eigenvector $\mathbf{q}_j$ ($j\neq i$) yields the fundamental relation
\begin{equation}\label{eq:proj_deriv}
    \mathbf{q}_j^H\mathbf{M}\mathbf{q}_i' =
    \frac{\mathbf{q}_j^H\mathbf{K}'\mathbf{q}_i}{\lambda_i - \lambda_j},
\end{equation}
which shows that the projection of the eigenvector derivative onto mode $j$ is
proportional to the inter‑modal coupling $\mathbf{q}_j^H\mathbf{K}'\mathbf{q}_i$
and inversely proportional to the eigenvalue gap $\lambda_i-\lambda_j$.  This
inverse‑gap dependence is the central reason why eigenvectors can change abruptly
when two eigenvalues approach each other.

Expanding $\mathbf{q}_i'$ in the $\mathbf{M}$-orthonormal eigenbasis and
enforcing the standard phase convention $\mathbf{q}_i^H\mathbf{M}\mathbf{q}_i'=0$
(to eliminate the self‑term) gives the explicit representation (detailed derivation provided in \cref{app:perturbation})
\begin{equation}\label{eq:qprime_expansion}
    \mathbf{q}_i' = \sum_{j\neq i}
    \frac{\mathbf{q}_j^H\mathbf{K}'\mathbf{q}_i}{\lambda_i - \lambda_j}\,\mathbf{q}_j .
\end{equation}
\Cref{eq:qprime_expansion} provides a unified, quantitative picture of modal
interaction.  In a veering region, a small eigenvalue gap combined with finite
coupling produces a large derivative, directly explaining the rapid eigenvector
rotation that breaks correlation‑based tracking.  Conversely, at a
symmetry‑protected crossing the coupling term vanishes identically (as required
by the Wigner–von Neumann non‑crossing rule); the derivative remains bounded and
the eigenvectors evolve smoothly.  These two scenarios represent the generic
outcomes of the non‑crossing rule for a single‑parameter Hermitian system, and
they motivate the adaptive, gap‑informed refinement strategy developed in
\cref{sec:adaptive}.

Finally, the $\mathbf{M}$-norm of the derivative,
\begin{equation}\label{eq:norm_deriv}
    \|\mathbf{q}_i'\|_{\mathbf{M}}^2 =
    \sum_{j\neq i} \frac{|\mathbf{q}_j^H\mathbf{K}'\mathbf{q}_i|^2}{(\lambda_i-\lambda_j)^2},
    \quad \|\mathbf{v}\|_{\mathbf{M}}^2 := \mathbf{v}^H\mathbf{M}\mathbf{v},
\end{equation}
provides a compact measure of the local sensitivity; it is well‑defined because
every denominator is non‑zero under the simplicity assumption.

\subsection{Challenges in SAFE mode tracking}
\label{sec:challenges}

While the MAC provides a simple tracking tool, its success depends sensitively
on the eigenvector evolution with $k$.  Mode crossing, veering, and degeneracy
can all cause misidentification.  Their essential interaction is captured by a
reduced two‑mode model.

\subsubsection{A two‑state local approximation}

Projecting \cref{eq:eigenproblem} onto the subspace spanned by two
near‑coincident eigenvectors $\mathbf{q}_i(k_0),\mathbf{q}_j(k_0)$ gives the
effective $2\times2$ Hermitian matrix
\begin{equation}\label{eq:2x2_H}
    \mathbf{H}(k)= \begin{pmatrix}
        a(k) & b(k)  \\
        b^{H}(k) & c(k)
    \end{pmatrix},
\end{equation}
with $a,c\in\mathbb{R}$ and $b\in\mathbb{C}$.  The diagonal entries are
$a(k)=\mathbf{q}_i^H(k_0)\mathbf{K}(k)\mathbf{q}_i(k_0)$,
$c(k)=\mathbf{q}_j^H(k_0)\mathbf{K}(k)\mathbf{q}_j(k_0)$, and the coupling
$b(k)=\mathbf{q}_i^H(k_0)\mathbf{K}(k)\mathbf{q}_j(k_0)$.  Because the
eigenvectors are $\mathbf{M}$‑orthonormal at $k_0$, $b(k_0)=0$; the first‑order
behaviour is governed by $b'(k_0)=\mathbf{q}_i^H\mathbf{K}'\mathbf{q}_j$, which
is precisely the coupling coefficient appearing in
\cref{eq:proj_deriv}.  The two eigenvalues of $\mathbf{H}$ are
\begin{equation}\label{eq:2x2_eigen}
    \lambda_{\pm}(k) = \frac{a+c}{2} \pm
    \frac{1}{2}\sqrt{|a-c|^2 + 4|b|^2}.
\end{equation}

\subsubsection{Four scenarios of mode interaction and their implications for MAC tracking}
\label{sec:four_scenarios}
\Cref{tab: modal_interactions} classifies the possible interactions in
single‑parameter Hermitian eigenproblems, linking each scenario to the eigengap,
coupling strength, and symmetry structure.  This classification directly informs
the tracking strategy.

\begin{table}[htbp]
\centering
\caption{Classification of modal interactions in single-parameter Hermitian eigenvalue problems arising from SAFE formulations, with physical interpretation and tracking implications. The local behavior near any interaction is captured by the two-state approximation \cref{eq:2x2_H,eq:2x2_eigen}, where the global eigenvalues $\lambda_i$, $\lambda_j$ are approximated by the local branches $\lambda_{\pm}$.}
\label{tab: modal_interactions}
\begin{tabular}{@{}p{2.0cm}p{2.0cm}p{1.5cm}p{4.7cm}p{3.8cm}@{}}
\toprule
\textbf{Scenario} & \textbf{Eigenvalue} & \textbf{Coupling} & \textbf{Physical mechanism} & \textbf{MAC tracking} \\
\midrule
Symmetry-protected crossing & $\lambda_i = \lambda_j$ at isolated $k$ & $b(k) \equiv 0$ & Orthogonal deformation mechanisms; no energy exchange despite eigenvalue coincidence & Easy: eigenvectors vary smoothly; bounded derivatives \\
\addlinespace
Avoided crossing (veering) & $\lambda_i \approx \lambda_j$ (minimum gap $2|b|>0$ ) & $b(k) > 0$ & Generic behavior in coupled modes; rapid exchange of physical character over narrow $k$-interval & Demands small $\Delta k$: derivative $\propto 1/\text{gap}$; misassignment risk high \\
\addlinespace
Crossing with non-zero coupling & $\lambda_i = \lambda_j$ & $b(k) \neq 0$ & Impossible in single-parameter Hermitian systems (Wigner--von Neumann rule) & --- \\
\addlinespace
Symmetry-protected degeneracy & $\lambda_i \equiv \lambda_j$ (all $k$) & $b(k) \equiv 0$ & Continuous symmetry (e.g., axisymmetry); eigenvectors rotate arbitrarily in invariant subspace & Pointwise MAC invalid; subspace tracking required \\
\bottomrule
\end{tabular}
\end{table}

\textbf{Scenario 1---Symmetry-protected crossing (SPC).}
When the two modes belong to different irreducible representations of a
symmetry respected by the waveguide, the off‑diagonal coupling vanishes
identically: $b(k)\equiv0$ for all $k$. The eigenvalues are simply $a(k)$ and
$c(k)$, and a true crossing occurs when $a=c$.  The eigenvectors do not mix and
evolve analytically, so the NCC~\cref{eq:NCC} holds even for
coarse steps---this is the benign situation analysed in \cref{sec:symmetry}.

\textbf{Scenario 2---Avoided crossing (veering).}
Without symmetry protection, $b(k)\neq0$ generically.  The eigenvalues repel,
reaching a minimum gap $2\min|b|>0$; over this narrow interval the eigenvectors
undergo a rapid but continuous exchange.  \Cref{eq:qprime_expansion} shows
that the eigenvector derivative scales as $1/\text{gap}$, explaining the
tracking difficulty.  A sufficiently small step $\Delta k$ always exists
(Proposition~1) and the adaptive algorithm of \cref{sec:adaptive} supplies it.

This is the generic behaviour in coupled systems.  In symmetric laminates,
mirror symmetry separates Lamb waves into $S$ and $A$ families that may cross,
but modes within the same family still veer.  In unsymmetric waveguides, all
modes are hybridised and \emph{every} interaction is an avoided crossing; no
pure modes exist, and veering is a ubiquitous feature of the spectrum.

\textbf{Scenario 3---Crossing with non‑zero coupling} is impossible in a
single‑parameter Hermitian system (Wigner--von Neumann rule) and is listed only
for completeness.

\textbf{Scenario 4---Symmetry-protected degeneracy (SPD).}
A continuous symmetry (e.g., axisymmetry) forces eigenvalues to coincide over
the entire wavenumber range ($a\equiv c$, $b\equiv0$).  The eigenvectors are not
unique but span a degenerate subspace in which any orthonormal basis is valid.
Pointwise MAC therefore fails: numerically computed basis vectors may rotate
arbitrarily, while the invariant subspace evolves smoothly.  Subspace tracking
(\cref{sec:spd}) is the physically correct approach.

These four scenarios capture the essential behaviour of mode interactions; the
following subsections provide the rigorous mathematical basis for an adaptive
algorithm that handles both symmetry‑protected crossings and veering robustly.

\subsection{Benign nature of symmetry-protected crossings}
\label{sec:symmetry}

When a waveguide possesses reflection symmetry, the governing operator
decomposes into block‑diagonal form with respect to the irreducible
representations of the symmetry group.  Modes from different symmetry classes
therefore evolve in decoupled subspaces, and the inter‑modal coupling vanishes
identically:
\begin{equation}\label{eq:symmetry_coupling}
    \mathbf{q}_j^H(k)\,\mathbf{K}'(k)\,\mathbf{q}_i(k) \equiv 0.
\end{equation}
This is a direct consequence of the Wigner--von Neumann non‑crossing rule
\cite{von_neumann_no_1929,hatton_noncrossing_1976,mead_noncrossing_1979}: a true
crossing $\lambda_i(k_0)=\lambda_j(k_0)$ can occur only between different
irreducible representations~\cite{gravenkamp_notes_2023}.

\subsubsection{Decoupled evolution and bounded eigenvector derivatives}
\label{sec:bounded_derivatives}

Because the coupling \cref{eq:symmetry_coupling} vanishes identically (not
just at the crossing point), the resonant term in the eigenvector derivative
expansion \cref{eq:qprime_expansion} that would diverge as
$1/(\lambda_i-\lambda_j)$ is absent.  No singular behaviour arises: the
eigenvector derivative remains bounded through the crossing, and each mode
family varies analytically within its decoupled symmetry block.

\subsubsection{Analyticity, Taylor expansion, and MAC tractability}
\label{sec:mac_tractability}

Standard perturbation theory~\cite{kato_perturbation_1995} guarantees that
eigenvectors within each block can be chosen analytic in $k$.  Expanding
$\mathbf{q}_i(k+\Delta k)$ in a Taylor series and using the phase convention
$\mathbf{q}_i^H\mathbf{M}\mathbf{q}_i'=0$ gives the self‑MAC decay
\begin{equation}
    \mathrm{MAC}[\mathbf{q}_i(k),\mathbf{q}_i(k+\Delta k)]
    = 1 - (\Delta k)^2\,\|\mathbf{q}_i'(k)\|_{\mathbf{M}}^2
      + O((\Delta k)^3),
    \label{eq:self_MAC_analytic}
\end{equation}
while the cross‑MAC with modes from other symmetry classes remains of
order $(\Delta k)^2$.  The Taylor expansion radius is finite and
independent of the eigenvalue gap; consequently, the
NCC~\cref{eq:NCC} can be satisfied with a step size
constrained by analyticity, not by the vanishing gap.

In contrast, at an avoided crossing the eigenvalue branches remain
analytic, but the eigenvector derivatives grow as $1/\text{gap}$.  The
effective Taylor radius therefore shrinks as the gap narrows, forcing
the step size to decrease correspondingly.  This fundamental
distinction---step size limited by a finite analyticity radius in
symmetry‑protected crossings, versus limited by the shrinking eigengap
in veering---explains why symmetry‑protected crossings require no
special refinement, while veering zones demand adaptive sampling.

\textbf{Remark.}  The situation is fundamentally different when the
symmetry forces eigenvalues to coincide over an entire interval ---
symmetry‑protected degeneracy.  Although the coupling also vanishes, the
eigenvectors are not unique within the degenerate subspace; this case is
analysed in \cref{sec:spd} and requires subspace tracking.

\subsection{Existence of critical step size in avoided-crossing regions}
\label{sec:existence}
We now address avoided crossings (veering), where the eigenvalue gap becomes
small while the coupling $\mathbf{q}_j^H\mathbf{K}'\mathbf{q}_i$ remains
non‑zero.  \Cref{eq:qprime_expansion} then implies large eigenvector
derivatives, so MAC‑based tracking requires a step $\Delta k$ sufficiently
small to resolve the rapid variation. The following proposition, adapted from
classical perturbation theory~\cite{kato_perturbation_1995}, shows that such a
step always exists under the analyticity and simplicity assumptions of
\cref{sec:perturbation}.

\medskip
\noindent\textbf{Proposition 1 (Existence of a critical step size).}
Consider the single‑parameter Hermitian eigenproblem
$[\mathbf{K}(k) - \lambda\mathbf{M}]\mathbf{q} = \mathbf{0}$ for
$k \in [k_{\min},k_{\max}] \subseteq \mathbb{R}$, with $\mathbf{K}(k)$
analytic and $\mathbf{M}$ constant, real and positive definite.  Assume the
eigenvalues are simple throughout the interval.  Then for every
$k \in [k_{\min},k_{\max})$ and every mode $i$, there exists a positive number
$\Delta k_{\max}(k)$, depending on the local eigengap and coupling strength,
such that for all $0 < \Delta k \le \Delta k_{\max}(k)$,
\begin{equation}\label{eq: theorem_inequality}
    \mathrm{MAC}[\mathbf{q}_i(k),\mathbf{q}_i(k+\Delta k)]
    > \mathrm{MAC}[\mathbf{q}_i(k),\mathbf{q}_j(k+\Delta k)], \quad \forall\,j\neq i.
\end{equation}
The column‑wise analogue follows by interval reversal. The derivation of this condition is provided in \cref{app:theorem1_proof}.

\textit{Remarks.}
\begin{enumerate}
    \item The bound $\Delta k_{\max}(k)$ shrinks where gaps are tiny,
    justifying adaptive refinement in veering zones.
    \item The proposition guarantees algorithmic consistency
    (the Hungarian algorithm returns the identity permutation), not physical
    correctness.  When the gap is comparable to numerical noise, the NCC may
    hold even though the physical modes have exchanged character; in such cases
    symmetry‑based constraints or further refinement are needed to ensure
    physical identification.
\end{enumerate}

\subsubsection{MAC separation and a posteriori error indicator}
\label{sec: MAC_separation}
Proposition~1 ensures that, for every mode $i$, the MAC separation
\begin{equation}\label{eq: MAC_sep}
\begin{aligned}
    D_i(k,\Delta k) &:=
    \min\!\bigl(D_i^{\mathrm{row}}(k,\Delta k),\,
              D_i^{\mathrm{col}}(k,\Delta k)\bigr),\\[2pt]
    D_i^{\mathrm{row}} &= \mathrm{MAC}[\mathbf{q}_i(k),\mathbf{q}_i(k+\Delta k)]
                         - \max_{j\neq i}\mathrm{MAC}[\mathbf{q}_i(k),\mathbf{q}_j(k+\Delta k)],\\
    D_i^{\mathrm{col}} &= \mathrm{MAC}[\mathbf{q}_i(k),\mathbf{q}_i(k+\Delta k)]
                         - \max_{j\neq i}\mathrm{MAC}[\mathbf{q}_j(k),\mathbf{q}_i(k+\Delta k)]
\end{aligned}
\end{equation}
is strictly positive for all $\Delta k \le \Delta k_{\max}(k)$.
Physically, $D_i \approx 1$ indicates adiabatic mode‑shape evolution, whereas
$D_i \ll 1$ signals that the sampling is too coarse to resolve the rapid
eigenvector rotation in a veering zone.  To obtain a single scalar measure for
the whole interval we define the \textbf{a posteriori error indicator}
\begin{equation}
    \varepsilon(k,\Delta k) = 1 - \min_i D_i(k,\Delta k),
    \label{eq:error_indicator}
\end{equation}
which lies in $[0,2]$.  A value $\varepsilon \approx 0$ implies that every mode
is tracked with high confidence; $\varepsilon \approx 1$ indicates that at
least one mode is at risk of misidentification.

\subsection{Symmetry‑protected degeneracy: subspace tracking}
\label{sec:spd}

When the waveguide possesses a continuous symmetry (e.g., axisymmetry),
eigenvalues may remain degenerate over the entire wavenumber range, as in the
$\sin(m\theta)$ and $\cos(m\theta)$ pair of a circular cylinder.  Although the
coupling $\mathbf{q}_j^H\mathbf{K}'\mathbf{q}_i$ vanishes identically, the
eigenvectors are not unique: any orthonormal basis of the invariant subspace is
valid.  Pointwise MAC tracking therefore becomes ill‑posed because numerically
computed basis vectors can rotate arbitrarily between successive steps even
while the subspace evolves smoothly.

\subsubsection{Subspace similarity measure}
Degeneracy is detected by eigenvalue clustering (within a tolerance) and
verified by checking that the projected coupling $\mathbf{Q}^H\mathbf{K}'\mathbf{Q}$
is numerically negligible.  Let $\mathbf{Q}(k)$ and $\mathbf{Q}(k+\Delta k)$ be
$n\times d$ matrices whose columns form $\mathbf{M}$-orthonormal bases of the
$d$-dimensional degenerate subspace at the two wavenumbers.  The subspace MAC,
a natural rotation‑invariant extension of the vector MAC, is defined as
\begin{equation}\label{eq:subspace_MAC}
    \mathrm{MAC}_{\mathrm{sub}}[\mathbf{Q}(k),\mathbf{Q}(k+\Delta k)]
    = \frac{1}{d}\,
      \bigl|\mathbf{Q}^H(k)\,\mathbf{M}\,\mathbf{Q}(k+\Delta k)\bigr|_F^2,
\end{equation}
where $|\cdot|_F$ is the Frobenius norm.  It takes values in $[0,1]$ and
reduces to the standard vector MAC when $d=1$.

\subsubsection{Unified treatment in the adaptive framework}
A degenerate subspace can be treated as a single entity in the adaptive
algorithm, exactly like a non‑degenerate mode.  The MAC separation for a
subspace $\mathcal{S}$ is
\begin{equation}\label{eq:subspace_sep}
\begin{aligned}
    D_{\mathcal{S}}(k,\Delta k) &=
    \min\!\bigl(D_{\mathcal{S}}^{\mathrm{row}},\,
              D_{\mathcal{S}}^{\mathrm{col}}\bigr),\\[2pt]
    D_{\mathcal{S}}^{\mathrm{row}} &=
    \mathrm{MAC}_{\mathrm{sub}}[\mathbf{Q}_{\mathcal{S}}(k),\mathbf{Q}_{\mathcal{S}}(k+\Delta k)]
    - \max_{\mathcal{T}\neq\mathcal{S}}
      \mathrm{MAC}_{\mathrm{sub}}[\mathbf{Q}_{\mathcal{S}}(k),\mathbf{Q}_{\mathcal{T}}(k+\Delta k)],\\[2pt]
    D_{\mathcal{S}}^{\mathrm{col}} &=
    \mathrm{MAC}_{\mathrm{sub}}[\mathbf{Q}_{\mathcal{S}}(k),\mathbf{Q}_{\mathcal{S}}(k+\Delta k)]
    - \max_{\mathcal{T}\neq\mathcal{S}}
      \mathrm{MAC}_{\mathrm{sub}}[\mathbf{Q}_{\mathcal{T}}(k),\mathbf{Q}_{\mathcal{S}}(k+\Delta k)],
\end{aligned}
\end{equation}
and the overall error indicator is $\varepsilon(k,\Delta k)=1-\min_{\mathcal{S}}D_{\mathcal{S}}(k,\Delta k)$,
precisely the same form as in \cref{eq:error_indicator}.  The core adaptive
logic (refinement if $\varepsilon>\bar{\varepsilon}$ and $\Delta k>\Delta k_{\min}$)
remains unchanged; only the initial clustering step to identify degenerate
subspaces is added.

\textbf{Physical interpretation.}
In an axisymmetric pipe, the $\sin(m\theta)$ and $\cos(m\theta)$ modes form a
degenerate pair with identical phase and group velocities.  A transducer of
finite azimuthal extent excites a fixed superposition --- a vector in the
invariant subspace --- rather than a pure basis function.  The subspace itself
evolves smoothly with $k$, whereas the numerical basis vectors within it may
rotate arbitrarily.  Subspace tracking therefore monitors the physically
observable quantity: the invariant subspace that a transducer can excite and a
receiver can detect.

\subsection{Two-level adaptive refinement: tracking consistency and physics-consistency}
\label{sec:adaptive}
The error indicator $\varepsilon$ guarantees numerical tracking consistency (the NCC), but as discussed in \cref{sec:NCC} it does not distinguish a genuine crossing from an unresolved veering. For waveguides with known structural symmetry, the Wigner--von Neumann non-crossing rule provides an additional constraint: modes of the same symmetry class cannot cross; any apparent crossing among them must be an unresolved veering.
Incorporating this physical knowledge directly into the refinement criterion
would create a circular dependency, because symmetry classes require resolved
eigenvectors, which in turn require refinement.  This epistemic hierarchy
motivates a two‑level strategy:

\begin{itemize}
    \item \textbf{Level~1 (tracking consistency).} The grid is refined using
    the purely numerical error indicator $\varepsilon$ until the NCC holds.
    No symmetry information is required, and the level applies to arbitrary
    waveguides.  Once converged, the eigenvectors are sufficiently smooth to
    permit reliable extraction of symmetry classes.
    \item \textbf{Level~2 (physics consistency).} Symmetry classes
    $\mathcal{C}(i)$ are extracted from the Level‑1 converged modes.  A
    symmetry penalty $\mathcal{S}_i$ is applied to any assignment that would
    imply a crossing between modes of the same class, and the
    symmetrically‑penalized error indicator $\tilde{\varepsilon}$ triggers
    additional refinement where numerical consistency masks a physically
    forbidden crossing.  Level~2 thus distinguishes symmetry‑allowed crossings
    (accepted) from unresolved avoided crossings (further refined).
\end{itemize}

\subsubsection{Level 1: NCC-driven refinement}
\label{sec:level1}
The grid is refined using $\varepsilon(k,\Delta k)$ until
$\varepsilon \le \bar{\varepsilon}$ or $\Delta k \le \Delta k_{\min}$.
This stage is purely numerical and requires no knowledge of mode symmetry
classes.

\subsubsection{Level 2: Symmetry-penalized physics refinement}
\label{sec:level2}
After Level~1 convergence, the symmetry class $\mathcal{C}(i)$ of each mode
is extracted (see \cref{sec:discussion_details}).  For each interval
$[k,k+\Delta k]$ and mode $i$, the symmetry penalty is
\begin{equation}
    \mathcal{S}_i(k,\Delta k) = \gamma\,
    \delta_{\mathcal{C}(i),\mathcal{C}(\sigma(i))}\,
    \mathbb{I}_{\text{cross}}(i,\sigma(i);k,\Delta k),
    \label{eq:sym_penalty_algorithm}
\end{equation}
where $\gamma \in [0,1]$ controls the constraint strength ($\gamma=1$ forbids
same‑class crossings), and the eigenvalue‑order reversal indicator is
\begin{equation}
    \mathbb{I}_{\text{cross}}(i;k,k+\Delta k) =
    \max_{\substack{j\neq i \\ \mathcal{C}(j)=\mathcal{C}(i)}}
    \mathbb{I}\!\left[
    \bigl(\lambda_i(k)-\lambda_j(k)\bigr)
    \bigl(\lambda_i(k+\Delta k)-\lambda_j(k+\Delta k)\bigr) < 0
    \right].
    \label{eq:sym_cross_indicator}
\end{equation}
The \textbf{symmetrically-penalized error indicator} (SPEI) is
\begin{equation}
    \tilde{\varepsilon}(k,\Delta k) =
    1 - \min_i \bigl[D_i(k,\Delta k) - \mathcal{S}_i(k,\Delta k)\bigr].
    \label{eq:SPEI}
\end{equation}
Since $\mathcal{S}_i \ge 0$, the SPEI can exceed $\bar{\varepsilon}$ even when
the NCC holds, triggering targeted secondary refinement.  Intervals reaching
$\Delta k_{\min}$ with $\tilde{\varepsilon} > \bar{\varepsilon}$ are flagged
as \emph{unresolved veering zones}.  The same logic applies to degenerate
subspaces, using the subspace error indicator and an analogous penalty.

\subsubsection{Algorithm description}
\label{sec:algorithm}

\begin{figure}[htb]
\centering
\includegraphics[width=0.75\columnwidth]{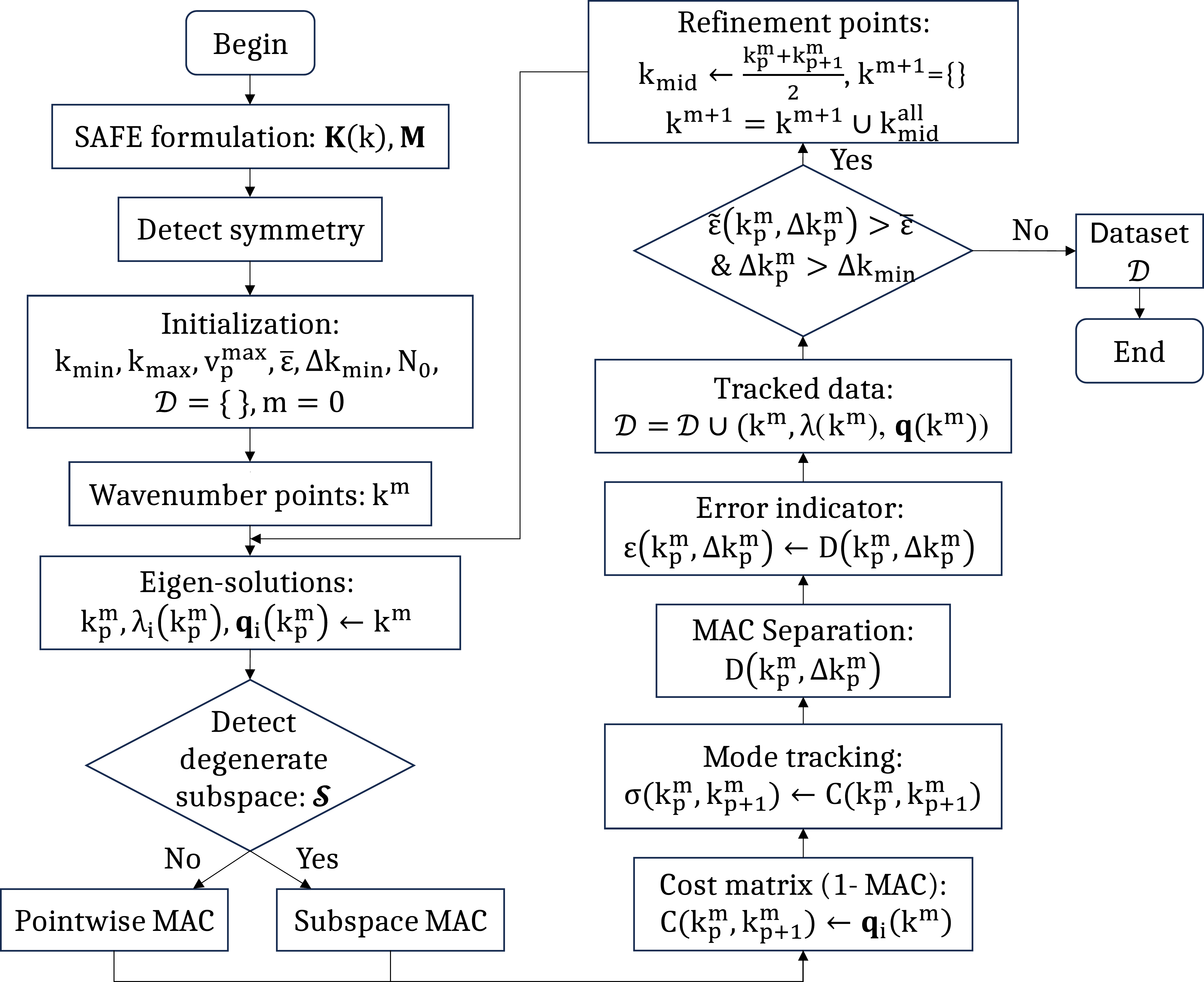}
\caption{Flowchart of the two‑level adaptive refinement algorithm. 
Level~1 ensures numerical tracking consistency (NCC); Level~2 enforces 
physical correctness via the symmetry-penalized error indicator (SPEI).}
\label{fig:adap_alg}
\end{figure}

The two‑level procedure is summarised in \cref{fig:adap_alg}.

\begin{enumerate}
    \item \textbf{Initialization.} Set $[k_{\min},k_{\max}]$,
    $\bar{\varepsilon}$, $\Delta k_{\min}$, $\gamma\in[0,1]$, and an initial
    uniform grid $\mathbf{k}^{(0)}$.
    \item \textbf{Solve eigenproblems.} At each $k_p^{(m)}$, compute all
    eigenpairs $(\lambda_i,\mathbf{q}_i)$.
    \item \textbf{Detect degenerate subspaces.} Cluster eigenvalues with
    negligible coupling; each cluster forms a degenerate subspace, triggering
    subspace MAC; otherwise use pointwise MAC.
    \item \textbf{Build cost matrix and assign.} For each interval, construct
    $\mathbf{C}$ with $C_{ij}=1-\mathrm{MAC}$ (or $1-\mathrm{MAC}_{\mathrm{sub}}$)
    and apply the Hungarian algorithm.
    \item \textbf{Level~1 refinement.} Compute $\varepsilon=1-\min_i D_i$.
    If $\varepsilon>\bar{\varepsilon}$ and $\Delta k_p>\Delta k_{\min}$, bisect
    the interval, update $\mathbf{k}^{(m+1)}$, and return to Step~2.
    \item \textbf{Level~1 convergence.} When no interval triggers refinement,
    proceed to Level~2.
    \item \textbf{Level~2: extract symmetry classes.} From the Level‑1
    converged eigenvectors, compute $\mathcal{C}(i)$ once; this classification
    remains fixed during subsequent iterations.  Evaluate
    $\tilde{\varepsilon}=1-\min_i(D_i-\mathcal{S}_i)$.
    \item \textbf{Level~2 refinement.} If $\tilde{\varepsilon}>\bar{\varepsilon}$
    and $\Delta k_p>\Delta k_{\min}$, bisect the interval and return to Step~2.
    Newly tracked modes inherit the symmetry class of their parent interval.
    \item \textbf{Termination.} Output the final grid; intervals reaching
    $\Delta k_{\min}$ without satisfying $\tilde{\varepsilon}\le\bar{\varepsilon}$
    are flagged as \emph{unresolved veering zones}.
\end{enumerate}

\subsubsection{Implementation details}
\label{sec:discussion_details}

\noindent\textbf{Symmetry class extraction and the role of $\gamma$.}
The symmetry penalty requires a mode symmetry class $\mathcal{C}(i)$ and the
coefficient $\gamma$.  For symmetric laminates, the parity of the
out‑of‑plane displacement classifies modes as S, A, or hybrid via
$\pi_i = \langle u_z(y,-z), u_z(y,z) \rangle / \|u_z\|^2$.  For axisymmetric pipes, the circumferential order $m$ is obtained from Fourier analysis.  Asymmetric waveguides assign
$\mathcal{C}(i)=\text{Trivial}$ to all modes.  In all cases we set
$\gamma=1$: for symmetric structures this prevents spurious crossings within
each symmetry class, while for asymmetric structures every mode belongs to the
same trivial class, so any eigenvalue‑order reversal is penalised and treated
as an unresolved veering.  The Kronecker delta
$\delta_{\mathcal{C}(i),\mathcal{C}(\sigma(i))}$ automatically restricts the
penalty to same‑class pairs.

\noindent\textbf{Convergence and tolerance.}
Bisection guarantees geometric step‑size reduction; Proposition~1 ensures a
finite local step exists under the simplicity assumption.  In all examples
tested, the process terminated well before reaching machine precision.  We
use $\bar{\varepsilon}=0.05$ throughout; reaching $\Delta k_{\min}$ without
Level~2 convergence signals a genuine veering whose gap is below numerical
resolution, not a failure.

\noindent\textbf{Computational overhead.}
Each refinement step solves eigenvalue problems only at newly inserted
wavenumbers; the total number of solves is proportional to the final grid
size, which is substantially smaller than a uniformly fine grid.  Symmetry‑class
extraction and SPEI evaluation add $O(N_{\text{modes}})$ operations, negligible
compared to the eigensolves.  The solves are independent and naturally parallel.

\section{Numerical results and discussions}
\label{sec:numerical}
The adaptive algorithm is validated on four examples covering symmetry‑protected
crossings, pervasive veering, arbitrary cross‑sections, and degeneracy: a
symmetric and an unsymmetric laminated plate, an L‑shaped aluminum bar, and a
steel pipe.  All are traction‑free.  Results are benchmarked against two
open‑source MATLAB tools that represent distinct tracking philosophies.
SAFEDC~\cite{liu_modelling_2024} uses a uniform wavenumber grid and
MAC‑based tracking; the Dispersion Calculator (DC)
v3.1~\cite{huber_dispersion_2024} employs a stiffness‑matrix formulation with
family‑based root‑finding.  DC first partitions modes into Lamb, shear
horizontal, etc., then solves $\det\mathbf{D}=0$ mode‑by‑mode with warm‑start
initialisation.  This works well for symmetric laminates but becomes
problematic when pure families do not exist, as in asymmetric configurations.
SAFEDC illustrates the limits of uniform sampling; DC provides an independent
validation of our MAC‑based tracking.  DC sweeps in frequency, the present
method sweeps in wavenumber, so the comparison focuses on tracking fidelity
rather than absolute branch counts.

All quantities are normalised by a characteristic length $a$ and shear wave
speed $c_T$ (\cref{tab:params_4examples}).  The adaptive strategy uses
$\bar{\varepsilon}=0.05$, $\Delta k_{\min}=10^{-4}$, an initial uniform grid
of $\Delta k=0.1$ ($10$ points per unit normalised wavenumber), and bisection
refinement when $\varepsilon>\bar{\varepsilon}$ (Level~1) or
$\tilde{\varepsilon}>\bar{\varepsilon}$ (Level~2).  Plate structures are
discretised with two Gauss–Lobatto–Legendre elements of order~5 per lamina
(identical to SAFEDC); the L‑bar and pipe use nine‑node quadrilaterals with
Gaussian quadrature.  The SAFE solver and adaptive algorithm are implemented
in Python; agreement with the independent MATLAB references confirms the
correctness of the implementation.
\begin{table}[ht]
\centering
\caption{Material properties, characteristic length $a$ and characteristic wave velocity $c_T$ for the four numerical examples.} \label{tab:params_4examples}
\small{
\begin{tabular}{llllll}
\toprule
Example & Material & Isotropy & Properties & Length $a$ & Velocity $c_T$ \\
\midrule
\begin{tabular}{@{}l@{}} Symmetric \\ laminate \end{tabular} 
& \multirow{2}{*}[-1ex]{M21T800}
& \multirow{2}{*}[-1ex]{\begin{tabular}{@{}l@{}} Transversely \\ isotropic \end{tabular}}
& \multirow{2}{*}{\begin{tabular}{@{}l@{}} $E_1=171$ GPa, $E_2=11.47$ GPa, \\ $G_{12}=4.83$ GPa, $\rho=1600$ $\mathrm{Kg/m^3}$, \\ $\nu_{12}=0.33$, $\nu_{23}=0.33$\end{tabular}} 
& \multirow{2}{*}[-1ex]{\begin{tabular}{@{}l@{}} Half laminate \\ thickness  \end{tabular}} 
& \multirow{2}{*}[-1ex]{3000 $\mathrm{m/s}$}
\\
\begin{tabular}{@{}l@{}} Unsymmetric \\ laminate \end{tabular} 
& 
& 
&  
&
& \\
\begin{tabular}{@{}l@{}} L‑shaped \\ bar \end{tabular} 
& Aluminium 
& Isotropic 
& \begin{tabular}{@{}l@{}} $E=70$ GPa, $\nu=0.33$, \\ $\rho=2700$ $\mathrm{Kg/m^3}$\end{tabular}  
& \begin{tabular}{@{}l@{}} Half short \\ leg length \end{tabular}  
& 3040\textsuperscript{1} $\mathrm{m/s}$ \\
\begin{tabular}{@{}l@{}} Steel \\ pipe \end{tabular} 
& Steel 
& Isotropic 
& \begin{tabular}{@{}l@{}} $E=207$ GPa, $\nu=0.3$, \\ $\rho=7850$ $\mathrm{Kg/m^3}$\end{tabular} 
& \begin{tabular}{@{}l@{}} Half radial \\ thickness \end{tabular} 
& 3000 $\mathrm{m/s}$ \\
\bottomrule
\end{tabular}} 
    \begin{tablenotes}
    \footnotesize
    \item[1] 1 The same $c_T$ is used as in Ref. \cite{maruyama_continuation_2025} for comparison purposes.
    \end{tablenotes}
\end{table}

\subsection{Symmetric laminated plate: symmetry-protected crossings and intra-family veering}
\begin{figure}[htb]
\centering
\includegraphics[width=1.03\columnwidth]{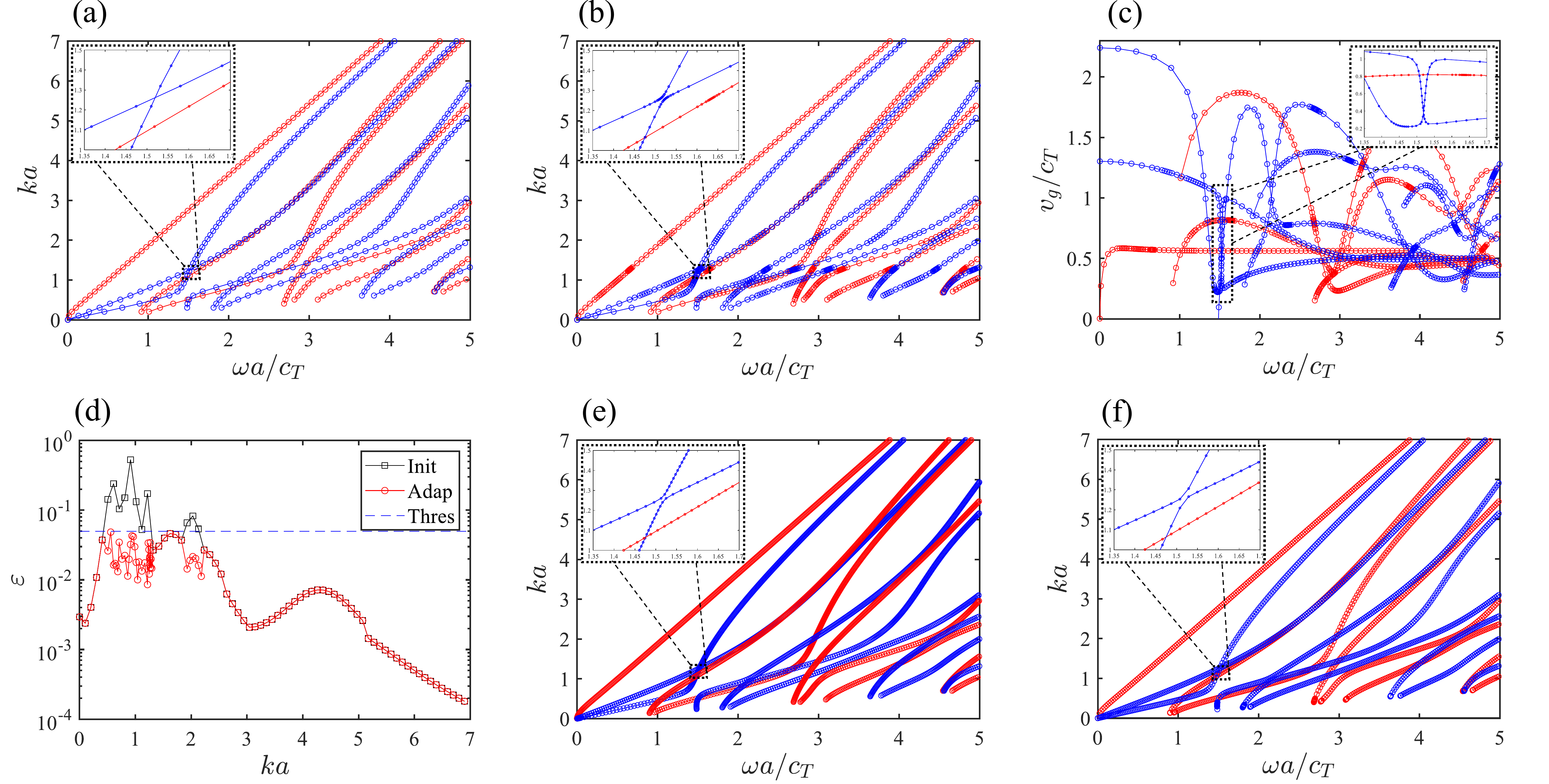}
\caption{Dispersion curves, group velocities, and tracking diagnostics for the symmetric laminated plate $[0,90,45,-45]_{2s}$. Symmetric (S) modes are shown in blue, antisymmetric (A) modes in red. (a)~Initial coarse sampling ($\Delta k = 0.1$, 70 points). (b)~Adaptive refinement result (98 points). (c)~Group velocity $v_g = \mathrm{d}\omega/\mathrm{d}k$ obtained from the adaptive grid. (d)~Error indicator $\varepsilon(k,\Delta k)$ on the initial and refined grids; the inset zooms the avoided crossing between two S-mode branches. (e)~SAFEDC with uniform fine sampling ($\Delta k = 0.02$, 350 points). (f)~DC results ($\Delta \omega = 0.0209$).}
\label{fig:exa1_sym_w_k_3solvers}
\end{figure}

The symmetric $[0,90,45,-45]_{2s}$ M21/T800 laminate decouples into symmetric (S)
and antisymmetric (A) families.  \Cref{fig:exa1_sym_w_k_3solvers} shows the
dispersion curves, group velocities, and error indicator.  DC confirms that
true crossings occur only between different symmetry classes; within the S
family the apparent intersection (inset of \cref{fig:exa1_sym_w_k_3solvers}(f))
is a classical avoided crossing.  This distinction carries direct practical
consequences: symmetric and antisymmetric modes involve orthogonal deformation
mechanisms—in‑plane stretching versus out‑of‑plane bending—so a transducer
that selectively excites S modes will traverse S–A crossings without modal
confusion.  The veering between two S modes, however, presents a genuine
hazard: two modes of the same symmetry class rapidly exchange their displacement
distributions over a narrow frequency band, and a fixed‑frequency inspection
protocol risks misidentifying which physical mode carries a defect signature.

The group‑velocity plot (\cref{fig:exa1_sym_w_k_3solvers}(c)) reinforces this
observation.  At $ka\approx1.5$ the two participating S branches swap $v_g$:
the branch that enters the interaction with a higher group velocity exits with
a lower value, and vice versa.  This behaviour follows directly from the
eigenvector rotation that characterises avoided crossings, and it provides an
independent, physically observable marker of veering.  For pulse‑echo
inspections that rely on time‑of‑flight for defect localisation, such a rapid
exchange implies that a slight shift in excitation frequency could route energy
along a fundamentally different propagation path, compromising the reliability
of source–defect distance estimates.

On the initial coarse grid ($\Delta k=0.1$, 70 points) all S–A crossings are
correctly tracked, but the two S‑mode branches are misidentified through the
veering region.  This failure illustrates a fundamental limitation of uniform
sampling: when the eigenvalue gap becomes small, the eigenvector variation over
a fixed step can exceed the local critical step size, so even a seemingly
reasonable global step cannot guarantee correct tracking.  SAFEDC with a
uniformly fine grid ($\Delta k=0.02$, 350 points) resolves the veering, but at
a five‑fold increase in cost and without any indication of where refinement is
actually needed or whether the assignments are reliable.

The adaptive algorithm automatically detects the veering region through a
pronounced peak of $\varepsilon$ on the coarse grid
(\cref{fig:exa1_sym_w_k_3solvers}(d)).  After local bisection the error indicator
falls below $\bar{\varepsilon}=0.05$, the S‑mode veering is correctly resolved
(\cref{fig:exa1_sym_w_k_3solvers}(b)), and the final grid contains only 98
points—substantially fewer than the uniform fine grid.  Crucially, the
post‑refinement error indicators provide quantitative confidence in the
tracking across the entire wavenumber range, a quality measure absent from
conventional tools.  This first example thus demonstrates that the adaptive
scheme handles symmetry‑protected crossings without special treatment,
automatically identifies inadequate step sizes in veering zones, and refines
only where needed.

\subsection{Unsymmetric laminated plate: pervasive veering}
\label{sec:unsym1}
The unsymmetric $[0,90,45,-45]_{4}$ M21/T800 laminate eliminates through‑thickness
symmetry, introducing tension–shear and bending–torsion coupling that hybridises
all modes.  Pure S, A, and SH families no longer exist, and by the Wigner–von
Neumann rule every modal interaction is an avoided crossing—veering is a
ubiquitous feature of the spectrum.  This makes the configuration a stringent
test for any tracking algorithm.

\begin{figure}[htb]
\centering
\includegraphics[width=1.03\columnwidth]{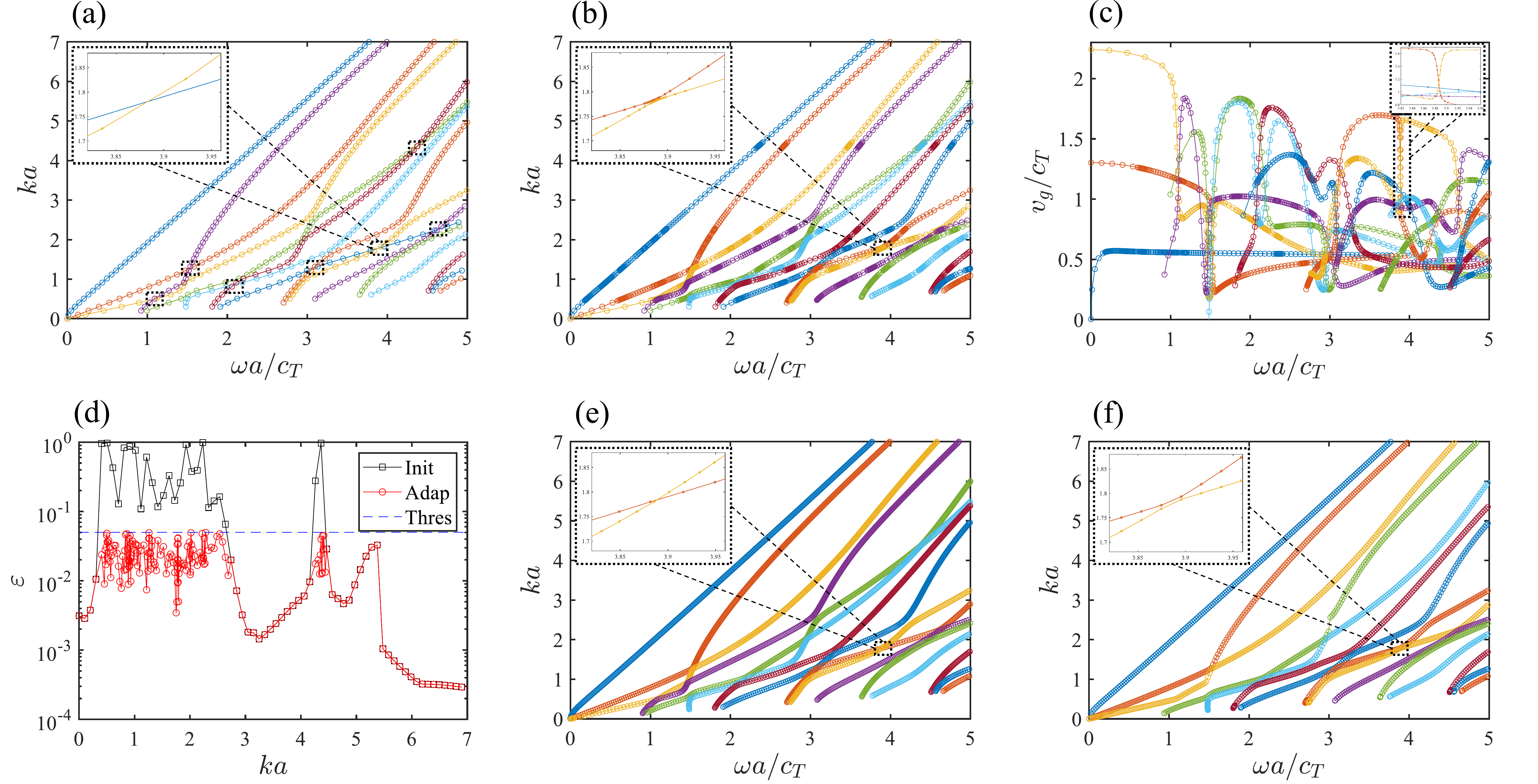}
\caption{Dispersion curves, group velocities, and tracking diagnostics for the unsymmetric laminated plate $[0,90,45,-45]_{4}$. (a)~Initial coarse sampling ($\Delta k = 0.1$, 70 points). (b)~Adaptive refinement result (185 points). (c)~Group velocity $v_g = \mathrm{d}\omega/\mathrm{d}k$ obtained from the adaptive grid. (d)~Error indicator $\varepsilon(k,\Delta k)$ on the initial and refined grids. (e)~SAFEDC with uniform fine sampling ($\Delta k = 0.02$, 350 points). (f)~DC results ($\Delta \omega = 0.0209$).}
\label{fig:exa2_unsym_w_k_3solvers}
\end{figure}

\Cref{fig:exa2_unsym_w_k_3solvers} collects the dispersion curves, group
velocities, and error indicator.  On the initial coarse grid ($\Delta k=0.1$,
70 points) multiple mode misidentifications occur (marked by black squares),
concentrated in veering zones.  SAFEDC with uniform $\Delta k=0.02$ (\cref{fig:exa2_unsym_w_k_3solvers}(e), 350 points)
reduces the errors to one but does not eliminate them entirely.  This residual
error underscores a fundamental limitation of uniform sampling: unless the step
size is chosen smaller than the minimum of $\Delta k_{\max}(k)$ over the entire
wavenumber range—a value that is unknown a priori and may be extremely small in
the most severe veering region—some misidentifications will persist.
Determining such a globally sufficient resolution is computationally prohibitive.

The DC results (\cref{fig:exa2_unsym_w_k_3solvers}(f)) confirm the complete
absence of true crossings.  The lack of symmetry classes also degrades DC's
family‑based root‑finding: without symmetry‑based pre‑partitioning to provide
initial guesses, the continuation corrector can become ill‑conditioned where
eigenvalue gaps are small, and even with a fine frequency step for mode
detection ($\Delta f=0.025$\,kHz), one mode is missed.  In contrast, the SAFE
formulation computes the complete eigenbasis at every wavenumber and therefore
captures all modes irrespective of coupling strength.

The group‑velocity plot (\cref{fig:exa2_unsym_w_k_3solvers}(c)) reveals that
arbitrary mode pairs exchange $v_g$ in veering zones, irrespective of nominal
polarisation—a signature of pervasive hybridisation.  For broadband pulse–echo
inspections, this implies that energy velocity can switch abruptly between
distinct wavefield families, compromising time‑of‑flight based defect
localisation unless the veering structure is fully resolved.

The error indicator on the coarse grid (\cref{fig:exa2_unsym_w_k_3solvers}(d))
exhibits pronounced peaks in $ka\in[0.3,2.6]$ and $[4.1,4.5]$, precisely where
misidentifications occur.  After recursive bisection, $\varepsilon$ falls
uniformly below $\bar{\varepsilon}=0.05$, and all veering zones are correctly
resolved with only 185 points—roughly half the uniform fine grid—and with
built‑in confidence measures.  The algorithm thus handles pervasive veering
without prior knowledge of mode families or critical step sizes, concentrating
points only where needed.

\subsection{L-shaped bar: arbitrary cross-section}
\label{sec:exa3}
\begin{figure}[htb]
\centering
\includegraphics[width=1.03\columnwidth]{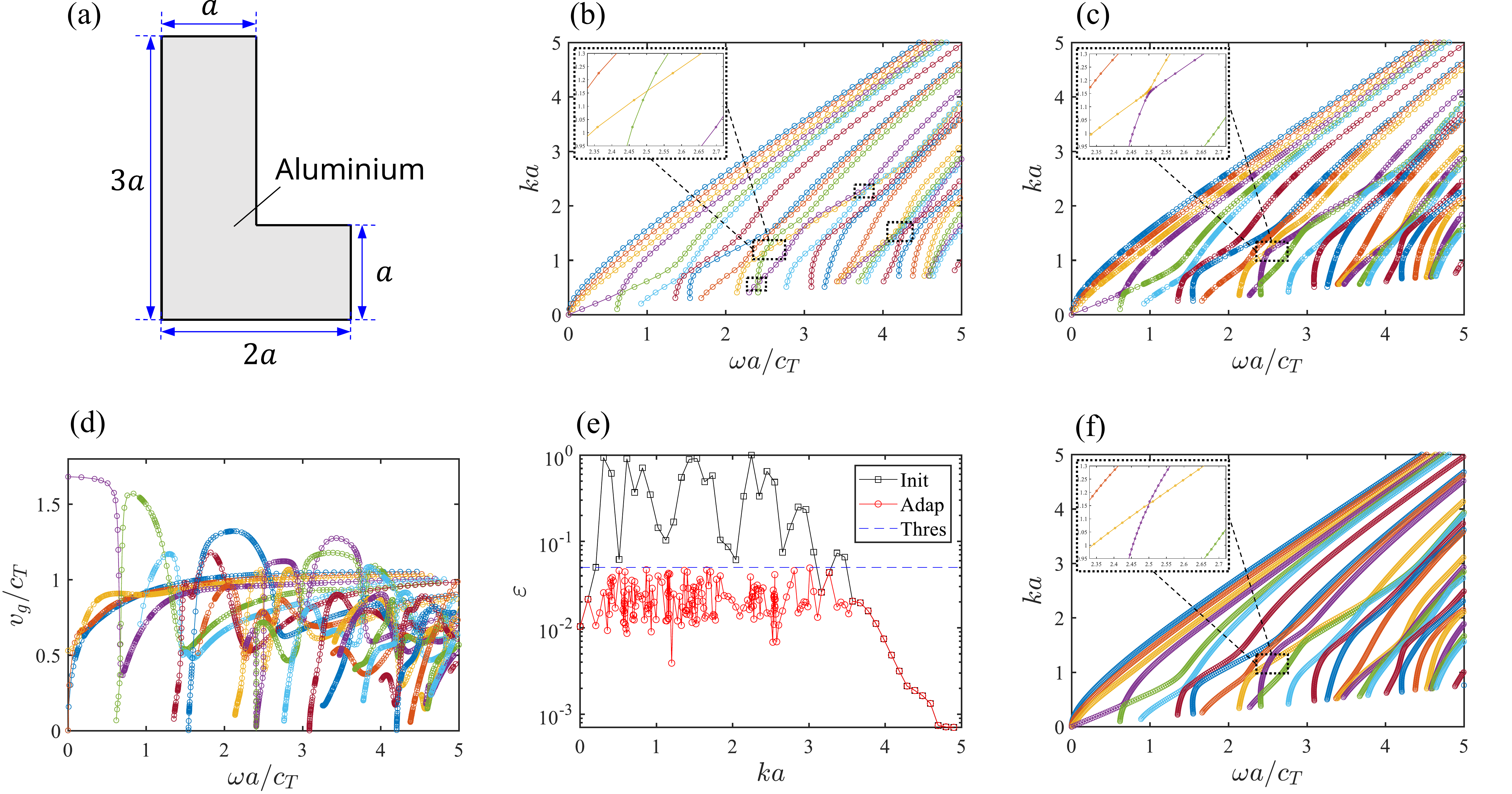}
\caption{Dispersion curves, group velocities, and tracking diagnostics for the L-shaped aluminum bar. (a)~Cross-section geometry. (b)~Initial coarse sampling ($\Delta k = 0.1$, 50 points). (c)~Adaptive refinement result (191 points). (d)~Group velocity $v_g = \mathrm{d}\omega/\mathrm{d}k$ obtained from the adaptive grid. (e)~Error indicator $\varepsilon(k,\Delta k)$ on the initial and refined grids. (f)~SAFEDC with uniform fine sampling ($\Delta k = 0.02$, 250 points). The reference NCM (Fig.~7 in Ref.~\cite{maruyama_continuation_2025}) confirms the absence of true crossings.}
\label{fig:exa3_Lbar_w_k_2solvers}
\end{figure}

The L‑shaped aluminum bar (cross‑section in \cref{fig:exa3_Lbar_w_k_2solvers}(a))
was recently benchmarked by Maruyama et al.~\cite{maruyama_continuation_2025},
whose numerical continuation method revealed pervasive mode veering and no true
crossings.  This is consistent with the Wigner–von Neumann rule, as the
cross‑section lacks the reflection symmetries that could give rise to
symmetry‑protected intersections.

\Cref{fig:exa3_Lbar_w_k_2solvers} shows the dispersion results.  The initial
coarse grid ($\Delta k=0.1$, 50 points) misidentifies several branches (marked
by black squares).  SAFEDC with uniform $\Delta k=0.02$ (250 points) reduces
but does not eliminate the errors; one erroneous branch persists.  This
residual error illustrates a practical limitation: without prior knowledge of
the most severe veering, even a globally fine grid can fail to resolve some
intervals completely.  The group‑velocity plot
(\cref{fig:exa3_Lbar_w_k_2solvers}(d)) displays the same pairwise $v_g$ exchange
observed in the unsymmetric plate, confirming the absence of symmetry
protection.

For structural health monitoring of industrial components—angles, channels, and
welded joints—geometric complexity precludes symmetry‑based mode classification.
In such structures, veering is the expected behaviour, and the inability to
predict where it occurs \emph{a priori} makes adaptive, confidence‑guided
sampling essential for reliable dispersion analysis.

The error indicator (\cref{fig:exa3_Lbar_w_k_2solvers}(e)) peaks in
$ka\in[0.2,3.4]$ on the coarse grid.  After adaptive bisection,
$\varepsilon$ drops below $\bar{\varepsilon}$ everywhere, all inaccurate
crossings are corrected to veering, and the final grid contains 191
points—fewer than the uniform fine grid and with guaranteed tracking
confidence.  This confirms that the adaptive strategy extends reliably to
arbitrary cross‑sections, automatically identifying and refining veering
regions without prior knowledge of the modal structure.

\subsection{Steel pipe: closed-section waveguide with symmetry-protected degeneracy}
\label{sec: pipe}
The steel pipe (inner radius $13a$, outer radius $15a$) is intentionally
modelled with a two‑dimensional Cartesian SAFE mesh, avoiding the customary
azimuthal Fourier decomposition.  This captures all circumferential orders
simultaneously, providing a stringent test of the algorithm's ability to
separate modes without prior knowledge of the underlying symmetry.  The
axisymmetry manifests as rotational invariance: for each circumferential
order $m>0$ the $\sin(m\theta)$ and $\cos(m\theta)$ modes form a degenerate
pair that shares the same eigenvalue for all $k$.  The coupling vanishes
identically, yet eigenvectors are unique only up to an arbitrary rotation
within the two‑dimensional invariant subspace (\cref{sec:spd}).  Pointwise MAC
tracking is therefore ill‑posed: numerically computed basis vectors can rotate
arbitrarily between adjacent $k$ steps, even though the subspace itself evolves
smoothly.

\begin{figure}[htb]
\centering
\includegraphics[width=1.03\columnwidth]{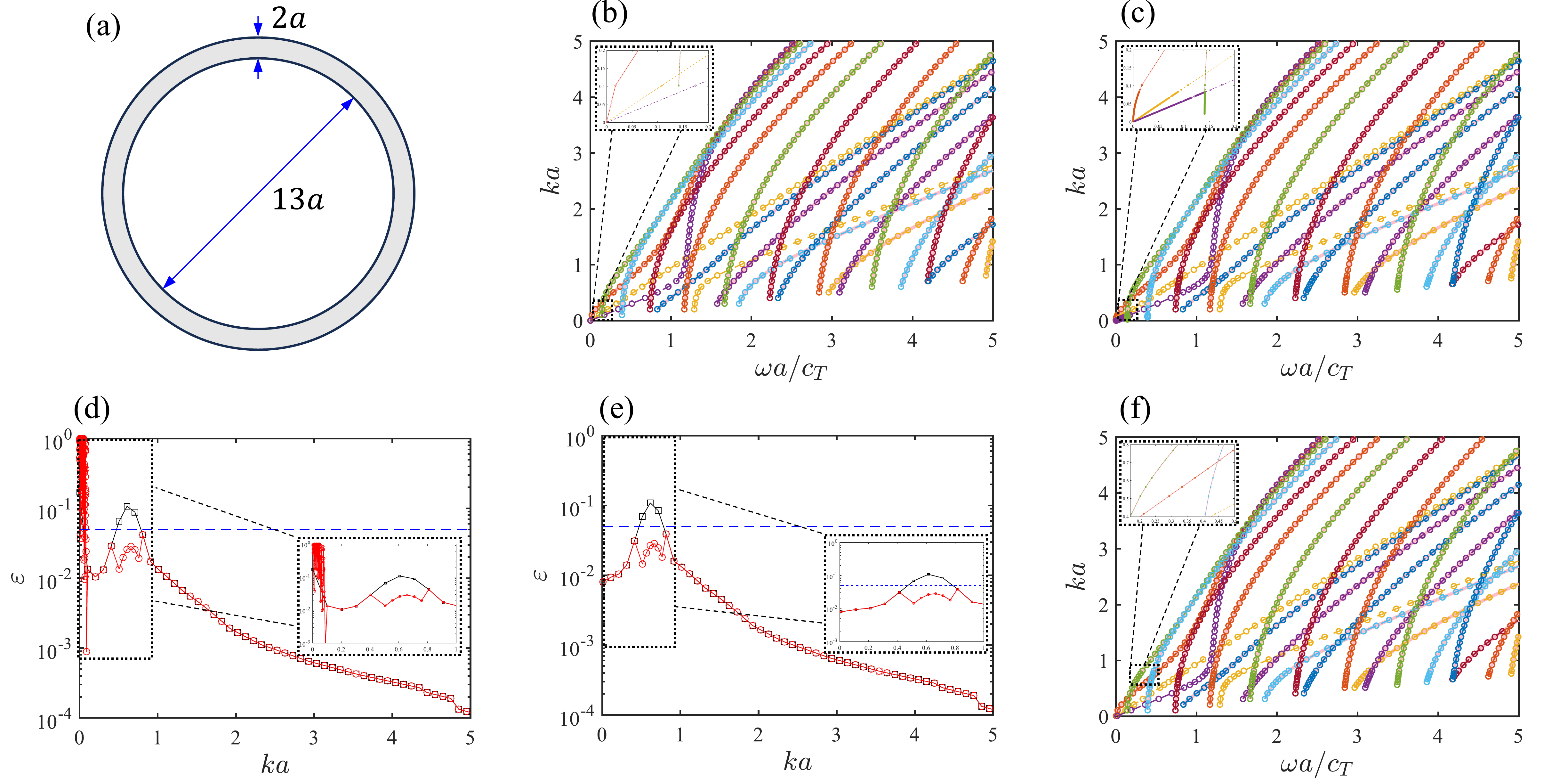}
\caption{Dispersion curves and tracking diagnostics for the steel pipe. (a)~Cross-section geometry. (b)~Initial coarse sampling ($\Delta k = 0.1$, 100 points). (c)~Adaptive refinement result (201 points). (d)~Error indicator $\varepsilon(k,\Delta k)$ with pointwise MAC tracking. (e)~Error indicator $\varepsilon(k,\Delta k)$ with subspace MAC tracking. (f)~SAFEDC with uniform fine sampling ($\Delta k = 0.04$, 250 points). For each circumferential order $m>0$, the $\sin(m\theta)$ and $\cos(m\theta)$ modes form a degenerate pair, producing coincident branches. All three dispersion strategies yield identical curves, validating the adaptive method's accuracy.}
\label{fig:exa4_annulus_w_k_2solvers}
\end{figure}

\Cref{fig:exa4_annulus_w_k_2solvers} presents the dispersion curves. In these dispersion curves, each SPD pair is distinguished by line style (e.g., solid versus dashed) for visual clarity, though the two components collapse onto the same dispersion branch. The three strategies (coarse, adaptive, uniform fine) yield identical curves,
confirming accuracy; multiple symmetry‑protected crossings between different
circumferential orders are visible and cause no difficulty.

The error indicators reveal a deeper story.  With pointwise MAC
(\cref{fig:exa4_annulus_w_k_2solvers}(d)) $\varepsilon$ cannot be reduced below
$\bar{\varepsilon}$ in the low‑$k$ regime, even after multiple refinement
iterations.  This does not indicate a physical discontinuity; rather, it
reflects the arbitrary rotation of the computed basis within the degenerate
subspace.  In the low‑$k$ regime the eigenproblem is dominated by the
cross‑sectional stiffness $\mathbf{K}_1$, and the numerical eigensolver has no
preferred direction inside the degenerate subspace.  Consequently,
infinitesimal perturbations—from round‑off, iterative tolerances, or mesh
discretisation—can rotate the computed basis arbitrarily between adjacent
$k$ steps, even though the invariant subspace itself evolves smoothly and the
degeneracy is exact.  As $k$ increases, the axial terms $ik\mathbf{K}_2$ and
$k^2\mathbf{K}_3$ modify the eigenproblem so that the degenerate subspace
deforms more rapidly.  This faster evolution imposes a tight differentiability
constraint on the eigenvectors: to maintain high similarity between successive
$k$‑steps, the computed basis is forced to follow the unique continuously
deforming pair that spans the subspace.  The pointwise MAC thus stabilises at
higher $k$, not because the degeneracy is lifted, but because the tracking
continuity condition becomes sufficiently restrictive.  This mechanism
clarifies why subspace tracking is indispensable at low $k$: there the
subspace evolves so slowly that no preferred basis is selected, whereas the
physical observable—the subspace itself—remains the same well‑defined invariant
under any rotation.

The adaptive algorithm detects degenerate eigenvalue clusters on the coarse
grid and switches to subspace MAC tracking.  Each degenerate pair is treated as
a single entity represented by its invariant subspace; the subspace MAC
provides a rotation‑invariant measure of similarity.  After adaptive refinement
guided by subspace MAC, the error indicator $\varepsilon$ falls below
$\bar{\varepsilon}$ across the entire wavenumber range
(\cref{fig:exa4_annulus_w_k_2solvers}(e)).

For pipeline inspection, this distinction is physically significant.  A
piezoelectric transducer of finite azimuthal extent excites a fixed
superposition within the degenerate subspace—neither pure $\sin(m\theta)$ nor
pure $\cos(m\theta)$, but a vector determined by the transducer's angular
position.  The experimentally observable quantity is the invariant subspace,
which evolves smoothly with $k$, while the numerical basis vectors returned by
the eigensolver may rotate freely.  Subspace tracking therefore monitors the
wavefields that an actual transducer can excite and a receiver can measure,
ensuring that the dispersion curves used for defect sizing correspond to
physically realisable wavefields.

\subsection{Computational efficiency}

A central motivation for adaptive sampling is the fundamental limitation of
uniform grids in resolving veering phenomena.  In a narrow veering zone the
eigengap becomes extremely small, requiring a step size smaller than the local
$\Delta k_{\max}(k)$ throughout the entire wavenumber range---yet
$\Delta k_{\max}(k)$ in the most severe veering region is unknown a priori
and may be orders of magnitude smaller than needed elsewhere.  The resulting
global refinement is computationally prohibitive.  \Cref{tab:efficiency}
shows that the adaptive final grids contain 98--201 points, compared to
250--350 for the uniform fine grids, and the disparity grows as veering zones
become narrower or more numerous.  Crucially, even these dense uniform grids
left residual misidentifications
(\cref{fig:exa2_unsym_w_k_3solvers}(e), \cref{fig:exa3_Lbar_w_k_2solvers}(f)),
whereas the adaptive grids guarantee numerical consistency by construction.
Because eigenvalue problems at different wavenumbers are independent, the
algorithm is naturally parallel; in the Python implementation, wall‑clock time
scaled roughly with the number of available cores.

\begin{table}[ht]
\centering
\caption{Comparison of grid sizes for the numerical examples.}
\label{tab:efficiency}
\small{
\begin{tabular}{llll}
\toprule
Example & Coarse initial  & Adaptive final & Uniform fine (SAFEDC) \\
\midrule
Symmetric laminate & 70 & 98 & 350 \\
Unsymmetric laminate & 70 & 185 & 350 \\
L‑shaped bar & 50 & 191 & 250 \\
Steel pipe & 100 & 201 & 250 \\
\bottomrule
\end{tabular}}
\end{table}

\subsection{Fundamental distinction: numerical matching versus physical consistency}
\label{sec:comparison}
The limitations observed in the preceding examples are not peculiar to SAFEDC
or DC, but reflect a common feature of existing tracking strategies: they treat
mode tracking as a purely numerical matching problem.  SAFEDC maximises
correlation between successive eigenvectors via MAC; DC traces classified mode
families by stepwise root‑finding with warm‑start initialisation; continuation
methods~\cite{maruyama_continuation_2025} follow parameterised solution curves
via predictor–corrector path‑following.  All three paradigms seek numerical
continuity---correlation, root convergence, or path proximity---without
verifying whether the resulting assignment respects the known physics of modal
interaction.  This suffices for well‑separated branches, but becomes inadequate
when the veering gap is small: the rapid eigenvector rotation over a narrow
$k$-interval can exceed the resolving power of a fixed step size, causing
misidentification.

Even when the step is small enough to satisfy the NCC and prevent label
jumping, a purely numerical tracker may still accept a numerically consistent
but physically impossible configuration.  If two modes of the same symmetry
class appear to cross on a grid that satisfies the NCC, the Wigner–von Neumann
non‑crossing rule dictates that this must be an unresolved veering rather than
a genuine crossing. A purely numerical criterion cannot detect this, because
it has no notion of symmetry classes. Only a physics‑informed criterion---such
as the symmetry penalty of our two‑level strategy---can flag such an apparent
crossing for further refinement. To our knowledge, the present work is the
first to incorporate the Wigner–von Neumann non‑crossing rule explicitly into
an automated mode‑tracking algorithm, thereby closing the gap between numerical
continuity and physical identification.

\subsection{Beyond exact symmetry: quasi-symmetric and damaged waveguides}
\label{sec:quasi_symmetric}

The preceding analysis focuses on waveguides possessing exact geometric and
material symmetry.  In practice, perfect symmetry is rarely realised.
Manufacturing tolerances, ply‑angle deviations, material ageing, environmental
variability, and impact damage introduce perturbations that weakly break the
underlying symmetry, causing real structures to occupy an intermediate regime
between ideal symmetry‑protected crossings and strongly asymmetric veering.

\subsubsection{From symmetry-protected crossings to modal hybridisation}

For an ideal symmetric waveguide, modes from different symmetry classes are
uncoupled, giving true crossings.  When symmetry is weakly broken, the coupling
between previously uncoupled families becomes non‑zero.  Any non‑zero coupling
removes the exact degeneracy and replaces the symmetry‑protected crossing with
an avoided crossing whose minimum eigengap scales with the induced coupling
strength.  The eigenvectors acquire admixtures of one another, with the degree
of mixing governed by the coupling strength.  Along the analytic eigenbranches,
modal identities are exchanged in the same manner as generic veering.  For weak
coupling the interaction region becomes asymptotically narrow, so the exchange
may be practically unobservable at finite numerical resolution.  Quasi‑symmetric
waveguides therefore represent the limiting case of mode veering, not a separate
class of interaction.

\subsubsection{Continuous transition and the role of adaptive sampling}

The transition from exact symmetry to strong asymmetry is continuous.  Exact
crossings occur only in the singular limit of vanishing coupling; any non‑zero
symmetry breaking produces an avoided crossing, while the observability of the
interaction depends on the relationship between the interaction width and the
numerical resolution.  The proposed adaptive algorithm handles all regimes
uniformly: it monitors local eigenvector variation through the MAC‑based error
indicator and refines the grid wherever the interaction becomes insufficiently
resolved, without requiring prior classification of the structural state.

\subsubsection{Implications for damaged waveguides}

Damage acts as a structural perturbation that may alter both eigenvalue
trajectories and modal coupling coefficients.  Symmetry‑preserving damage may
shift crossing locations but leave the protected nature intact;
symmetry‑breaking damage introduces coupling that transforms protected
crossings into avoided crossings, potentially too narrow to resolve
numerically.  The same perturbation‑theoretic framework therefore applies to
damaged structures: regions of strong eigenvector variation are refined
automatically, regardless of whether the underlying cause is geometric
asymmetry, manufacturing variability, or damage‑induced perturbations.  For
structures without a priori known symmetry classes, every mode belongs to the
same trivial symmetry class, and with $\gamma=1$ any apparent crossing is
treated as an unresolved veering.  An interval that reaches $\Delta k_{\min}$
without Level~2 convergence is honestly flagged as a \emph{below‑resolution
veering zone} rather than silently misidentified.  This transparent diagnostic
distinguishes between ``resolved and physically consistent'' and ``unresolved
but physically constrained,'' extending the framework from idealised waveguides
to the imperfect structures encountered in NDE and SHM practice.

\subsection{A guiding principle}
For conservative elastic waveguides the governing operator is self‑adjoint (the stiffness matrix is Hermitian or real symmetric), and the mode‑interaction phenomenology is dictated by structural symmetry alone.  Without symmetry, there are no symmetry‑protected crossings or degeneracies; within a mode family, veering is the generic behaviour, and its rapid eigenvector rotation scales inversely with the local eigenvalue gap, demanding adaptive resolution. With symmetry, crossings and degeneracies are protected by vanishing coupling and can be tracked reliably with coarse steps.  The resulting hierarchical logic—first verify that the system is Hermitian, then assess its symmetry, and finally interpret the local eigengap—provides the engineer with a systematic way to judge whether a computed dispersion diagram is physically consistent and where tracking demands finer resolution.

\section{Conclusion and future work}
\label{sec:conclusion}
This paper has applied classical perturbation theory to the SAFE eigenproblem,
establishing a quantitative relationship between eigenvector sensitivity and
eigenvalue gap. The analysis targets conservative elastic waveguides---systems
whose governing operator is self‑adjoint, rendering the stiffness matrix
Hermitian for lossless, traction‑free configurations.  This central insight
explains why correlation‑based tracking fails in veering regions while
remaining reliable at symmetry‑protected crossings.  Building on it, we
formulated a numerical consistency condition (NCC), proved the existence of a
critical step size below which the NCC holds, and developed a two‑level
adaptive refinement strategy that separates numerical tracking consistency from
symmetry‑based physical correctness.  The main conclusions are:

\begin{enumerate}
    \item \textbf{The eigengap controls eigenvector sensitivity.} The
    eigenvector derivative scales inversely with the eigenvalue gap, providing
    a unified interpretation of mode veering and the theoretical foundation for
    adaptive refinement.

    \item \textbf{Symmetry‑protected crossings are analytically benign.}
    Symmetry‑induced block‑diagonalisation forces the modal coupling to vanish
    identically; eigenvector derivatives remain bounded, so coarse‑step MAC
    tracking is reliable without special treatment.

    \item \textbf{A two‑level adaptive strategy separates numerical tracking
    from physical identification.} Level~1 refines the grid using the NCC until
    label continuity is guaranteed; Level~2 applies a symmetry penalty to
    correct numerically consistent but physically impossible assignments.
    Across four benchmark waveguides, this strategy resolves all veering zones
    with substantially fewer points than uniform fine sampling, while providing
    a built‑in error indicator as a quantitative confidence measure.

    \item \textbf{Degenerate subspaces require rotation‑invariant tracking.}
    In axisymmetric waveguides, pointwise MAC fails because eigenvectors rotate
    arbitrarily within each degenerate pair.  The subspace MAC treats the
    invariant subspace as a single entity, restores robust tracking, and
    reveals a wavenumber‑dependent stability that further justifies adaptive
    sampling.

    \item \textbf{Physical consistency is essential for reliable mode
    identification.} Traditional tools treat mode tracking as a purely
    numerical matching problem.  Incorporating the Wigner–von Neumann
    non‑crossing rule provides the missing physical constraint, particularly
    in asymmetric structures where all modes are coupled and every apparent
    crossing must be treated as a veering event.
\end{enumerate}

The entire framework is built on the premise that the system is Hermitian---a
condition satisfied by lossless, traction‑free waveguides, for which the
Wigner–von Neumann non‑crossing rule, the perturbation analysis, and the
existence guarantees hold rigorously.  The method is therefore not directly
applicable to non‑Hermitian systems such as viscoelastic or fluid‑loaded
waveguides, for which these theoretical foundations must be revisited.
Additional limitations include the heuristic tolerance $\bar{\varepsilon}$,
which may require tuning for different mode densities, the minimum step
$\Delta k_{\min}$ that can halt refinement of extremely narrow veering zones,
and the reliance on accurate eigenvector computation near cut‑off frequencies.

Future work should extend the framework to non‑Hermitian eigenproblems,
develop data‑driven tolerance selection based on local spectral properties, and
combine wavenumber adaptivity with $hp$‑refinement of the cross‑sectional mesh
for large‑scale heterogeneous structures.  A documented open‑source
implementation (TopoDisper) is available at
\url{https://github.com/dongxiao96/TopoDisper}, with a graphical user
interface under development to facilitate adoption by practitioners.

In summary, by systematically applying perturbation theory to the SAFE
eigenproblem, this work provides a physically interpretable and computationally
efficient solution to the long‑standing problem of mode‑tracking failure in
veering regions, bridging the gap between theoretical guarantees and
engineering practice.

\appendix
\section{Appendix: Detailed derivation of eigenvector perturbation analysis}
\label[appendix]{app:perturbation}

This appendix provides the complete step‑by‑step derivation of the eigenvector
derivative formulas presented in \cref{sec:perturbation}.

\subsection*{Differentiation of the eigenproblem}
Differentiating $\mathbf{K}\mathbf{q}_i = \lambda_i\mathbf{M}\mathbf{q}_i$ with respect to $k$ gives
\begin{equation}\label{eq:app:eigen_diff}
    \mathbf{K}'\mathbf{q}_i + \mathbf{K}\mathbf{q}_i' = \lambda_i'\mathbf{M}\mathbf{q}_i + \lambda_i\mathbf{M}\mathbf{q}_i',
\end{equation}
where primes denote $k$‑derivatives.  $\mathbf{K}'$ is well defined because $\mathbf{K}(k)$ is analytic (polynomial) in $k$.

\subsection*{Projection onto other modes}
For $j \neq i$, left‑multiply \cref{eq:app:eigen_diff} by $\mathbf{q}_j^H$:
\begin{equation}\label{eq:app:eigen_diff_lm}
    \mathbf{q}_j^H\mathbf{K}'\mathbf{q}_i + \mathbf{q}_j^H\mathbf{K}\mathbf{q}_i' = \lambda_i'\mathbf{q}_j^H\mathbf{M}\mathbf{q}_i + \lambda_i\mathbf{q}_j^H\mathbf{M}\mathbf{q}_i'.
\end{equation}
Using the adjoint relation $\mathbf{q}_j^H\mathbf{K} = \lambda_j\mathbf{q}_j^H\mathbf{M}$ (Hermitian $\mathbf{K}$) and $\mathbf{M}$‑orthogonality $\mathbf{q}_j^H\mathbf{M}\mathbf{q}_i = 0$ for $j\neq i$, we obtain
\begin{equation}
    \mathbf{q}_j^H \mathbf{K}' \mathbf{q}_i + \lambda_j \mathbf{q}_j^H \mathbf{M} \mathbf{q}_i' = \lambda_i \mathbf{q}_j^H \mathbf{M}\mathbf{q}_i',
\end{equation}
hence the fundamental projection formula
\begin{equation}\label{eq:app:proj_deriv}
    \mathbf{q}_j^H\mathbf{M}\mathbf{q}_i' = \frac{\mathbf{q}_j^H\mathbf{K}'\mathbf{q}_i}{\lambda_i - \lambda_j}, \quad j \neq i,
\end{equation}
which is \cref{eq:proj_deriv} in the main text.

\subsection*{Expansion of the eigenvector derivative}
Expand $\mathbf{q}_i'$ in the $\mathbf{M}$-orthonormal basis $\{\mathbf{q}_j\}_{j=1}^n$:
\begin{equation}\label{eq:app:expand_raw}
    \mathbf{q}_i' = \sum_{j=1}^n (\mathbf{q}_j^H\mathbf{M}\mathbf{q}_i') \mathbf{q}_j.
\end{equation}
Differentiating the normalisation $\mathbf{q}_i^H\mathbf{M}\mathbf{q}_i = 1$ gives
\begin{equation}
    \mathbf{q}_i^H\mathbf{M}\mathbf{q}_i' + (\mathbf{q}_i^H\mathbf{M}\mathbf{q}_i')^{H} = 2\,\mathrm{Re}(\mathbf{q}_i^H\mathbf{M}\mathbf{q}_i') = 0,
\end{equation}
so the real part vanishes; the imaginary part reflects the arbitrary $k$-dependent phase of the eigenvector.

\subsection*{Phase convention and boundedness}
We remove the self‑term $j=i$ in \cref{eq:app:expand_raw} by fixing the phase.  Multiplying each $\mathbf{q}_i$ by a $k$-dependent phase factor $\mathrm{e}^{\mathrm{i}\theta_i(k)}$ (which leaves eigenvalues and MAC unchanged) allows us to enforce
\begin{equation}\label{eq:app:phase_condition}
    \mathbf{q}_i^H\mathbf{M}\mathbf{q}_i' = 0.
\end{equation}
Such a choice preserves analyticity~\cite{kato_perturbation_1995}.  Inserting \cref{eq:app:proj_deriv} into \cref{eq:app:expand_raw} and using \cref{eq:app:phase_condition} yields
\begin{equation}\label{eq:app:qprime_expansion}
    \mathbf{q}_i' = \sum_{j \neq i} \frac{\mathbf{q}_j^H\mathbf{K}'\mathbf{q}_i}{\lambda_i - \lambda_j}\,\mathbf{q}_j,
\end{equation}
which is \cref{eq:qprime_expansion} in the main text.  The $\mathbf{M}$-norm follows directly:
\begin{equation}\label{eq:app:norm_deriv}
    \|\mathbf{q}_i'\|_{\mathbf{M}}^2 = \sum_{j \neq i} \frac{|\mathbf{q}_j^H\mathbf{K}'\mathbf{q}_i|^2}{(\lambda_i - \lambda_j)^2},
\end{equation}
reproducing \cref{eq:norm_deriv}.  Under the simplicity assumption every denominator is non‑zero, so the derivative is bounded away from true degeneracies.

\section{Proof of Proposition~1 (existence of a critical step size)}
\label[appendix]{app:theorem1_proof}

We prove that for each $k$ and mode $i$ there exists $\Delta k_{\max}(k)>0$ such that
\[
\mathrm{MAC}[\mathbf{q}_i(k),\mathbf{q}_i(k+\Delta k)] > \mathrm{MAC}[\mathbf{q}_i(k),\mathbf{q}_j(k+\Delta k)] \quad \forall\,j\neq i
\]
holds for all $0<\Delta k\le\Delta k_{\max}(k)$.  The proof follows from the analyticity and simplicity assumptions stated in \cref{sec:perturbation}.

\subsection*{Taylor expansion and MAC asymptotics}
Because eigenvectors can be chosen analytic~\cite{kato_perturbation_1995}, we expand
\begin{equation}\label{eq:app:taylor}
    \mathbf{q}_i(k+\Delta k) = \mathbf{q}_i(k) + \Delta k\,\mathbf{q}_i'(k) + \tfrac{1}{2}(\Delta k)^2\mathbf{q}_i''(k) + O((\Delta k)^3).
\end{equation}
Using $\mathbf{M}$-orthonormality and the phase convention $\mathbf{q}_i^H\mathbf{M}\mathbf{q}_i'=0$, the inner products become
\[
\mathbf{q}_i^{H}(k)\mathbf{M}\mathbf{q}_i(k+\Delta k) = 1 - \tfrac{1}{2}(\Delta k)^2\|\mathbf{q}_i'(k)\|_{\mathbf{M}}^2 + O((\Delta k)^3),
\]
and, for $j\neq i$,
\[
\mathbf{q}_i^{H}(k)\mathbf{M}\mathbf{q}_j(k+\Delta k) = \Delta k\,c_{ij}(k) + O((\Delta k)^2), \quad c_{ij}(k) = \frac{\mathbf{q}_i^H\mathbf{K}'\mathbf{q}_j}{\lambda_j - \lambda_i}.
\]
Hence the MAC values are
\begin{equation}\label{eq:app:self_MAC}
    \mathrm{MAC}[\mathbf{q}_i(k),\mathbf{q}_i(k+\Delta k)] = 1 - (\Delta k)^2\|\mathbf{q}_i'(k)\|_{\mathbf{M}}^2 + O((\Delta k)^3),
\end{equation}
\begin{equation}\label{eq:app:cross_MAC}
    \mathrm{MAC}[\mathbf{q}_i(k),\mathbf{q}_j(k+\Delta k)] = (\Delta k)^2|c_{ij}(k)|^2 + O((\Delta k)^3).
\end{equation}

\subsection*{Local estimates and construction of $\Delta k_{\max}(k)$}
From \cref{eq:app:self_MAC} and \cref{eq:app:cross_MAC} there exist $\delta_0>0$ and constants $C_1, C_2>0$ such that for all $0<\Delta k\le\delta_0$ and $j\neq i$,
\begin{equation}\label{eq:app:estimates}
    \mathrm{MAC}[\mathbf{q}_i(k),\mathbf{q}_i(k+\Delta k)] \ge 1 - C_1(\Delta k)^2, \qquad
    \mathrm{MAC}[\mathbf{q}_i(k),\mathbf{q}_j(k+\Delta k)] \le C_2(\Delta k)^2.
\end{equation}
Define
\begin{equation}\label{eq:app:Delta_k_max}
    \Delta k_{\max}(k) := \min\!\left(\delta_0,\;\frac{1}{\sqrt{2(C_1+C_2)}}\right).
\end{equation}
For $0<\Delta k\le\Delta k_{\max}(k)$ we have $(\Delta k)^2 \le 1/[2(C_1+C_2)]$, which implies $1 - C_1(\Delta k)^2 > C_2(\Delta k)^2$.  Combining with \cref{eq:app:estimates} yields, for every $j\neq i$,
\[
\mathrm{MAC}[\mathbf{q}_i(k),\mathbf{q}_i(k+\Delta k)] \ge 1 - C_1(\Delta k)^2 > C_2(\Delta k)^2 \ge \mathrm{MAC}[\mathbf{q}_i(k),\mathbf{q}_j(k+\Delta k)],
\]
which is exactly the row‑wise NCC \cref{eq:NCC}.  The column‑wise condition follows by symmetry.  

\section*{CRediT authorship contribution statement}
\textbf{Dong Xiao}: Conceptualization, Methodology, Software, Data curation, Formal analysis, Writing - Original draft preparation, Writing - Review Editing, Visualization. \textbf{Zahra Sharif-Khodaei}: Supervision, Writing - Review Editing. \textbf{M. H. Aliabadi}: Supervision, Writing - Review Editing.

\section*{Declaration of competing interest}
The authors declare that they have no known competing financial interests or personal relationships that could have appeared to influence the work reported in this paper.

\section*{Acknowledgements}
The first author acknowledges the financial support from the K. C. Wong Postdoctoral Fellowship, funded by the K. C. Wong Education Foundation.

\section*{Data availability}
The source code and data supporting this study will be made publicly available upon publication at \url{https://github.com/dongxiao96/TopoDisper}. This Python-based tool implements the proposed adaptive wavenumber sampling algorithm, and the repository contains Jupyter notebooks that reproduce all numerical examples presented in this paper. A graphical user interface (GUI) is under development to further lower the barrier for users.


\setstretch{1.0}
\small{

}

\end{document}